\documentclass[global]{svjour}
\usepackage{amssymb}
\usepackage{amsfonts}
\usepackage{amsmath}

\input{tcilatex}

\begin{document}

\title{Fuzzy equilibrium existence for Bayesian abstract fuzzy economies and
applications to random quasi-variational inequalities with random fuzzy
mappings}
\author{Monica Patriche}
\institute{University of Bucharest 
\email{monica.patriche@yahoo.com}%
}
\mail{\\
University of Bucharest, Faculty of Mathematics and Computer Science, 14
Academiei Street, 010014 Bucharest, Romania}
\maketitle

{\small \ }\textbf{Abstract }In this paper, we introduce a Bayesian abstract
fuzzy economy model and we prove the Bayesian fuzzy equilibrium existence.
As applications, we prove the existence of the solutions for two types of
random quasi-variational inequalities with random fuzzy mappings and we also
obtain random fixed point theorems.


AMS\ Subject Classification: 58E35, 47H10, 91B50, 91A44.

\section{\textbf{\ INTRODUCTION}}

The study of fuzzy games has begun with the paper written by Kim and Lee in
1998 [16]. This type of games is a generalization of classical abstract
economies. For an overview of results concerning this topic, the reader is
referred to [24]. Though, the existence of random fuzzy equilibrium has not
been studied by now. We introduce the new model of Bayesian abstract fuzzy
economy and explore the existence of the Bayesian fuzzy equilibrium. Our
model is characterized by a private information set, an action (strategy)
fuzzy mapping, a random fuzzy constraint one and a random fuzzy preference
mapping. The Bayesian fuzzy equilibrium concept is an extension of the
deterministic equilibrium. We generalize the former deterministic models
introduced by Debreu [8], Shafer and Sonnenschein [25], Yannelis and
Prabhakar [29] or Patriche [24] and search for applications.

Since Fichera and Stampacchia introduced the variational inequalities (in
1960s), this domain has been extensively studied. For recent results we
refer the reader to [1]-[4], [6],[7], [12], [17], [20]-[22], [26], [28] and
the bibliography therein. Noor and Elsanousi [19] introduced the notion of a
random variational inequality. Existence of solutions of the random
variational inequality and random quasi-variational inequality problems has
been proved, for instance, in [12], [13], [18], [27], [33].

In this paper, we first define the model of the Bayesian abstract fuzzy
economy and we prove a theorem of Bayesian fuzzy equilibrium existence. Then
we apply it in order to prove the existence of solutions for two types of
random quasi-variational inequalities with random fuzzy mappings. We
generalize some results obtained by Yuan in [27]. As a consequence, we
obtain random fixed point theorems.

The paper is oragnized as follows. In the next section, some notational and
terminological conventions are given. We also present, for the reader's
convenience, some results on Bochner integration. In Section 3, the model of
differential information abstract fuzzy economy is introduced and the main
result is stated. Section 4 contains existence results for solutions of
random quasi-variational inequalities with random fuzzy mappings.

\section{\textbf{\protect\smallskip NOTATION AND DEFINITION}}

Throughout this paper, we shall use the following notation:

1. $%
\mathbb{R}
_{++}$ denots the set of strictly positive reals. co$D$ denotes the convex
hull of the set $D$. $\overline{co}D$ denotes the closed convex hull of the
set $D$. $2^{D}$ denotes the set of all non-empty subsets of the set $D$. If 
$D\subset Y$, where $Y$ is a topological space, cl$D$ denotes the closure of 
$D$.\smallskip

For the reader's convenience, we review a few basic definitions and results
from continuity and measurability of correspondences and Bochner integrable
functions.

Let $Z$ and $Y$ be sets.

\begin{definition}
The \textit{graph} of the correspondence $P:Z\rightarrow 2^{Y}$ is the set $%
G_{P}=\{(z,y)\in Z\times Y:y\in P(z)\}$.
\end{definition}

Let $Z$, $Y$ be topological spaces and $P:Z\rightarrow 2^{Y}$ be a
correspondence.

\QTP{Body Math}
1. $P$ is said to be \textit{upper semicontinuous} if for each $z\in Z$ and
each open set $V$ in $Y$ with $P(z)\subset V$, there exists an open
neighborhood $U$ of $z$ in $Z$ such that $P(y)\subset V$ for each $y\in U$.

2. $P$ is said to be \textit{lower semicontinuous} if for each $z\in Z$ and
each open set $V$ in $Y$ with $P(z)\cap V\neq \emptyset $, there exists an
open neighborhood $U$ of $z$ in $Z$ such that $P(y)\cap V\neq \emptyset $
for each $y\in U$.\medskip

\begin{lemma}
(see [32]). \textit{Let }$Z$\textit{\ and }$Y$\textit{\ be two topological
spaces and let }$D$\textit{\ be an open subset of }$Z.$\textit{\ Suppose }$%
P_{1}:Z\rightarrow 2^{Y}$\textit{\ , }$P_{2}:Z\rightarrow 2^{Y}$\textit{\
are upper semicontinuous correspondences such that }$P_{2}(z)\subset
P_{1}(z) $\textit{\ for all }$z\in D.$\textit{\ Then the correspondence }$%
P:Z\rightarrow 2^{Y}$\textit{\ defined by}
\end{lemma}

\begin{center}
$P\mathit{(z)=}\left\{ 
\begin{array}{c}
P_{1}(z)\text{, \ \ \ \ \ \ \ if }z\notin D\text{, } \\ 
P_{2}(z)\text{, \ \ \ \ \ \ \ \ \ \ if }z\in D%
\end{array}%
\right. $
\end{center}

\textit{is also upper semicontinuous.\medskip }

\begin{definition}
Let $Y$ be a metric space and $Y^{\prime }$ be its dual. $P:Y\rightarrow
2^{Y^{\prime }}$ is said to be\textit{\ monotone }if Re$\langle
u-v,y-x\rangle \geq 0$\textit{\ }for all $u\in P(y)$\ and $v\in P(x)$ and $%
x,y\in Y.\medskip $
\end{definition}

Let now $(\Omega ,$ $\tciFourier $, $\mu )$ be a complete, finite measure
space, and $Y$ be a topological space.

1. The correspondence $P:\Omega \rightarrow 2^{Y}$ is said to have a \textit{%
measurable graph} if $G_{P}\in \tciFourier \otimes \beta (Y)$, where $\beta
(Y)$ denotes the Borel $\sigma $-algebra on $Y$ and $\otimes $ denotes the
product $\sigma $-algebra.

2. The correspondence $T:\Omega \rightarrow 2^{Y}$ is said to be \textit{%
lower measurable} if for every open subset $V$ of $Y$, the set $\{\omega \in
\Omega $ $:$ $T(\omega )\cap V\neq \emptyset $\} is an element of $%
\tciFourier $.

Recall (see Debreu [9], p. 359) that if $T:\Omega \rightarrow 2^{Y}$ has a
measurable graph, then $T$ is lower measurable. Furthermore, if $T(\cdot )$
is closed valued and lower measurable then $T:\Omega \rightarrow 2^{Y}$ has
a measurable graph.

\begin{lemma}
\textbf{(}see\textbf{\ [}15\textbf{]).} Let $P_{n}:\Omega \rightarrow 2^{Y}$%
, $n=1,2...$be a sequence of correspondences with measurable graphs. Then
the correspondences $\cup _{n}P_{n}$, $\cap _{n}P_{n}$ and $Y\setminus P_{n}$
have measurable graphs.\medskip
\end{lemma}

Let $(\Omega ,\tciFourier $, $\mu )$ be a measure space and $Y$ be a Banach
space.

It is known (see [15], Theorem 2, p.45) that, if $x:\Omega \rightarrow Y$ is
a $\mu $-measurable function, then $x$ is Bochner integrable if only if $%
\underset{\Omega }{\int }\Vert x(\omega )\Vert d\mu (\omega )<\infty $.

It is denoted by $L_{1}(\mu ,Y)$ the space of equivalence classes of $Y$%
-valued Bochner integrable functions $x:\Omega \rightarrow Y$ normed by $%
\parallel x\parallel =\underset{\Omega }{\int }\Vert x(\omega )\Vert d\mu
(\omega )$. Also it is known (see [9], p.50) that \textit{\ }$L_{1}(\mu ,Y)$%
\textit{\ }is a Banach space.

\begin{definition}
The correspondence $P:\Omega \rightarrow 2^{Y}$ is said to be \textit{%
integrably bounded} if there exists a map $h\in L_{1}(\mu ,R)$ such that $%
sup\{\parallel x\parallel $ $:$ $x\in P(\omega )\}\leq h(\omega )$ $\mu -a.e$%
.
\end{definition}

We denote by $S_{P}^{1}$ the set of all selections of the correspondence $%
P:\Omega \rightarrow 2^{Y}$ that belong to the space $L_{1}$($\mu ,Y)$, i.e.

$S_{P}^{1}=\{x\in L_{1}(\mu ,Y):$ $x(\omega )\in P(\omega )$ $\mu $%
-a.e.\}.\medskip

Further, we will see the conditions under which $S_{P}^{1}$ is nonempty and
weakly compact in $L_{1}(\mu ,Y)$. Aumann measurable selection theorem (see
Appendix) and Diestel's Theorem (see Appendix) are necessary.\medskip

Let $\tciFourier (Y)$ be a collection of all fuzzy sets over $Y.$

\begin{definition}
A mapping $P:\Omega \rightarrow \tciFourier (Y)$ is called a fuzzy mapping.
If $P$ is a fuzzy mapping from $\Omega ,$ $P(\omega )$ is a fuzzy set on $Y$
and $P(\omega )(y)$ is the membership function of $y$ in $P(\omega ).$
\end{definition}

Let $A\in \tciFourier (Y),$ $a\in \lbrack 0,1],$ then the set $%
(A)_{a}=\{y\in Y:A(y)\geq a\}$ is called an $a-$cut set of fuzzy set $A.$

\begin{definition}
A fuzzy mapping $P:\Omega \rightarrow \tciFourier (Y)$ is said to be
measurable if for any given $a\in \lbrack 0,1],$ $(P(\cdot ))_{a}:\Omega
\rightarrow 2^{Y}$ is a measurable set-valued mapping.
\end{definition}

\begin{definition}
We say that a fuzzy mapping $P:\Omega \rightarrow \tciFourier (Y)$ is said
to have a measurable graph if for any given $a\in \lbrack 0,1],$ the
set-valued mapping $(P(\cdot ))_{a}:\Omega \rightarrow 2^{Y}$ has a
measurable graph.
\end{definition}

\begin{definition}
A fuzzy mapping $P:\Omega \times X\rightarrow \tciFourier (Y)$ is called a
random fuzzy mapping if for any given $x\in X,$ $P(\cdot ,x):\Omega
\rightarrow \tciFourier (Y)$ is a measurable fuzzy mapping.
\end{definition}

\section{BAYESIAN FUZZY EQUILIBRIUM EXISTENCE FOR BAYESIAN ABSTRACT FUZZY
ECONOMIES}

\subsection{THE\ MODEL OF A BAYESIAN ABSTRACT FUZZY ECONOMY}

We now define the next model of the Bayesian abstract fuzzy economy which
generalizes the model in [23].

Let $(\Omega $\textit{, }$\tciFourier ,\mu )$ be a complete finite measure
space, where $\Omega $ denotes the set of states of nature of the world and
the $\sigma -$algebra $\tciFourier $, denotes the set of events. Let $Y$
denote the strategy or com$\func{mod}$ity space, where $Y$ is a separable
Banach space.\medskip

Let $I$ be a countable or uncountable set (the set of agents). Let $%
X_{i}:\Omega \rightarrow \mathcal{F}(Y)$ a fuzzy mapping and $z\in (0,1].$

Let $L_{X_{i}}=\{x_{i}\in S_{(X_{i}(\cdot ))_{z}}:x_{i}$ is $\tciFourier
_{i} $ -measurable $\}.$ Denote by $L_{X}=\tprod\limits_{i\in I}L_{X_{i}}$
and by $L_{X_{-i}}$ the set $\tprod\limits_{j\neq i}L_{X_{j}}.$ An element $%
x_{i}$ of $L_{X_{i}}$ is called a strategy for agent $i$. The typical
element of $L_{X_{i}}$ is denoted by $\widetilde{x}_{i}$ and that of $%
X_{i}(\omega )$ by $x_{i}(\omega )$ (or $x_{i}$).

\begin{definition}
A \textit{general Bayesian abstract fuzzy economy} is a family $G=\{(\Omega
,\tciFourier ,\mu ),$ $(X_{i},\tciFourier
_{i},A_{i},P_{i},a_{i},b_{i},z_{i})_{i\in I},\}$, where
\end{definition}

(1) $X_{i}:\Omega \rightarrow \mathcal{F}(Y)$ is the \textit{action
(strategy) fuzzy mapping} of agent $i$,

(2) $\tciFourier _{i}$ is a sub $\sigma -$algebra of $\tciFourier $ which
denotes the \textit{private information of agent i},

(3) for each $\omega \in \Omega ,$ $A_{i}(\omega ,$\textperiodcentered $%
):L_{X}\rightarrow \mathcal{F}(Y)$ is the \textit{random fuzzy constraint
mapping of agent }$i;$

(4) for each $\omega \in \Omega ,$ $P_{i}(\omega ,$\textperiodcentered $%
):L_{X}\rightarrow \mathcal{F}(Y)$ is the \textit{random fuzzy preference
mapping of agent }$i;$

(5) $z\in (0,1]$ is such that for all $(\omega ,x)\in \Omega \times L_{X},$ $%
(A_{i}(\omega ,\widetilde{x}))_{a_{i}(\widetilde{x})}\subset (X_{i}(\omega
))_{z}$ and $(P_{i}(\omega ,\widetilde{x}))_{p_{i}(\widetilde{x})}\subset
X_{i}(\omega ))_{z};$

(6) $a:S_{X}^{1}\rightarrow (0,1]$ is a random fuzzy constraint function and 
$p:S_{X}^{1}\rightarrow (0,1]$ is a random fuzzy preference function.$%
\medskip $

\begin{definition}
A \textit{Bayesian fuzzy equilibrium for }$G$ is a strategy profile $%
\widetilde{x}^{\ast }\in L_{X}$ such that for all $i\in I$,
\end{definition}

(j) $\widetilde{x}_{i}^{\ast }(\omega )\in $cl$(A_{i}(\omega ,\widetilde{x}%
^{\ast }))_{a_{i}(\widetilde{x}^{\ast })}$ $\mu -a.e.$

(jj) $(A_{i}(\omega ,\widetilde{x}^{\ast }))_{a_{i}(\widetilde{x}^{\ast
})}\cap (P_{i}(\omega ,\widetilde{x}^{\ast }))_{p_{i}(\widetilde{x}^{\ast
})}=\emptyset $ $\mu -a.e.\medskip $

\begin{remark}
Now we assume that for each $i\in I,$ $X_{i}$ is a compact convex nonempty
subset of $Y$ and for each $\omega \in \Omega ,$ we set $(X_{i}(\omega
))_{z}=X_{i}$. Then we obtain the deterministic classical model of
Yannelis-Prabhakar in [29] for an abstract economy with any set of players.
\end{remark}

\begin{remark}
The interpretation of the preference fuzzy mapping $P_{i}$ is that $y_{i}\in
(P_{i}(\omega ,\widetilde{x}))_{p_{i}(\widetilde{x})}$ means that at the
state $\omega $ of the nature, agent $i$ strictly prefers $y_{i}$ to $%
\widetilde{x}_{i}(\omega )$ if the given strategy of other agents is fixed.
The preference do not need to be reprezentable by utility functions.
However, it will be assumed that $\widetilde{x}_{i}(\omega )\notin
(P_{i}(\omega ,\widetilde{x}))_{p_{i}(\widetilde{x})}$ $\mu -a.e.\medskip $
\end{remark}

\subsection{EXISTENCE OF THE BAYESIAN FUZZY EQUILIBRIUM}

This is our first theorem. The constraint and preference correspondences
derived from the constraint and preference fuzzy mappings verify the
assumptions of measurable graph and weakly open lower sections. Our results
is a generalization of Theorem 3 in [23].\medskip

\begin{theorem}
\textit{Let }$I$\textit{\ be a countable or uncounatble set. Let the family} 
$G=\{(\Omega ,\tciFourier ,\mu ),$ $(X_{i},\tciFourier
_{i},A_{i},P_{i},a_{i},b_{i},z_{i})_{i\in I},\}$ \textit{be a general
Bayesian abstract economy satisfying A.1)-A.4). Then there exists a Bayesian
fuzzy equilibrium for }$G$\textit{.}
\end{theorem}

\textit{For each }$i\in I:$

\textit{A.1)}

\textit{\ \ \ \ \ (a) }$X_{i}:\Omega \rightarrow \mathcal{F}(Y)$\textit{\ is
such that }$\omega \rightarrow X_{i}(\omega )_{z}:\Omega \rightarrow 2^{Y}$ 
\textit{is} \textit{a nonempty, convex, weakly compact-valued and integrably
bounded correspondence.}

\textit{\ \ \ \ \ (b) }$X_{i}:\Omega \rightarrow \mathcal{F}(Y)$\textit{\ is
such that }$\omega \rightarrow (X_{i}(\omega ))_{z}:\Omega \rightarrow 2^{Y}$
\textit{is }$\tciFourier _{i}-$\textit{lower measurable;}

\textit{\ A.2)}

\textit{\ \ \ \ \ (a) For each }$(\omega ,\widetilde{x})\in \Omega \times
L_{X},$\textit{\ }$(A_{i}(\omega ,\widetilde{x}))_{a_{i}(\widetilde{x})}$%
\textit{\ is\ convex and has a non-empty interior in the relative norm
topology of }$(X_{i}(\omega ))_{z}.$\textit{\medskip }

\textit{\ \ \ \ \ (b)\ the correspondence }$(\omega ,\widetilde{x}%
)\rightarrow (A_{i}(\omega ,\widetilde{x}))_{a_{i}(\widetilde{x})}:\Omega
\times L_{X}\rightarrow 2^{Y}$\textit{\ has measurable graph i.e. }$%
\{(\omega ,\widetilde{x},y)\in \Omega \times L_{X}\times Y:y\in
(A_{i}(\omega ,\widetilde{x}))_{a_{i}(\widetilde{x})}\}\in \tciFourier
\otimes $\textit{\ss }$_{w}(L_{X})\otimes $\textit{\ss }$(Y)$\textit{\ where 
\ss }$_{w}(L_{X})$\textit{\ is the Borel }$\sigma -$\textit{algebra for the
weak topology on }$L_{X}$\textit{\ and \ss }$(Y)$\textit{\ is the Borel }$%
\sigma -$\textit{algebra for the norm topology on }$Y$\textit{.}

\textit{\ \ \ \ \ (c)\ the correspondence }$(\omega ,\widetilde{x}%
)\rightarrow (A_{i}(\omega ,\widetilde{x}))_{a_{i}(\widetilde{x})}$\textit{\
has weakly open lower sections, i.e., for each }$\omega \in \Omega $\textit{%
\ and for each }$y\in Y,$\textit{\ the set }$((A_{i}(\omega ,\widetilde{x}%
)_{a_{i}(\widetilde{x})})^{-1}(\omega ,y)=\{\widetilde{x}\in L_{X}:y\in
(A_{i}(\omega ,\widetilde{x}))_{a_{i}(\widetilde{x})}\}\}$\textit{\ is
weakly open in }$L_{X};$

\textit{\ \ \ \ \ (d) For each }$\omega \in \Omega ,$\textit{\ }$\widetilde{x%
}\rightarrow $cl$(A_{i}(\omega ,\widetilde{x}))_{a_{i}(\widetilde{x}%
)}:L_{X}\rightarrow 2^{Y}$\textit{\ is upper semicontinuous in the sense
that the set }$\{\widetilde{x}\in L_{X}\mathit{:}$cl$(A_{i}(\omega ,%
\widetilde{x}))_{a_{i}(\widetilde{x})}\subset V$ \textit{is weakly open in }$%
L_{X}$\textit{\ for every norm open subset }$V$\textit{\ of} $Y$.

\textit{A.3)}

\textit{\ \ \ \ \ (a) the correspondence }$(\omega ,\widetilde{x}%
)\rightarrow (P_{i}(\omega ,\widetilde{x}))_{p_{i}(\widetilde{x})}:\Omega
\times L_{X}\rightarrow 2^{Y}$\textit{\ has nonempty open convex values such
that }$(P_{i}(\omega ,\widetilde{x}))_{p_{i}(\widetilde{x})}\subset
(X(\omega ))_{z}$\textit{\ for each }$(\omega ,\widetilde{x})\in \Omega
\times L_{X}.$

\textit{\ \ \ \ (b)\ the correspondence} $(\omega ,\widetilde{x})\rightarrow
(P_{i}(\omega ,\widetilde{x}))_{p_{i}(\widetilde{x})}:\Omega \times
L_{X}\rightarrow 2^{Y}$ \textit{has measurable graph \ }

\textit{(c) the correspondence} $(\omega ,\widetilde{x})\rightarrow
(P_{i}(\omega ,\widetilde{x}))_{p_{i}(\widetilde{x})}:\Omega \times
L_{X}\rightarrow 2^{Y}$\textit{\ has weakly open lower sections, i.e., for
each }$\omega \in \Omega $\textit{\ and for each }$y\in Y,$\textit{\ the set 
}$((P_{i}(\omega ,\widetilde{x}))_{p_{i}(\widetilde{x})})^{-1}(\omega ,y)=\{%
\widetilde{x}\in L_{X}:y\in (P_{i}(\omega ,\widetilde{x}))_{p_{i}(\widetilde{%
x})})\}$\textit{\ is weakly open in }$L_{X};$

\textit{A.4)}

\textit{\ \ \ \ (a) For each }$\widetilde{x}_{i}\in L_{X_{i}},$\textit{\ for
each }$\omega \in \Omega ,$\textit{\ }$\widetilde{x}_{i}(\omega )\notin
(A_{i}(\omega ,\widetilde{x}))_{a_{i}(\widetilde{x})}\cap (P_{i}(\omega ,%
\widetilde{x}))_{p_{i}(\widetilde{x})}.\medskip $

\textit{Proof.} For each $i\in I,$ define $\Phi _{i}:\Omega \times
L_{X}\rightarrow 2^{Y}$ by $\Phi _{i}(\omega ,\widetilde{x})=A_{i}(\omega ,%
\widetilde{x}))_{a_{i}(\widetilde{x})}\cap P_{i}(\omega ,\widetilde{x}%
))_{p_{i}(\widetilde{x})}.$ We prove first that $L_{X}$ is a non-empty,
convex, weakly compact subset in $L_{1}(\mu ,Y).$

Since $(\Omega $\textit{, }$\tciFourier ,\mu )$ is a complete finite measure
space, $Y$ is a separable Banach space and $X_{i}:\Omega \rightarrow 2^{Y}$
has measurable graph, by Aumann's selection theorem (see Appendex) it
follows that there exists a $\tciFourier _{i}$-measurable function $%
f_{i}:\Omega \rightarrow Y$ such that $f_{i}(\omega )\in X_{i}(\omega )$ $%
\mu -a.e$. Since $X_{i}$ is integrably bounded, we have that $f_{i}\in
L_{1}(\mu ,Y)$, hence $L_{X_{i}}$ is non-empty and $L_{X}=\tprod\limits_{i%
\in I}L_{X_{i}}$ is non-empty. Obviously $L_{X_{i}}$ is convex and $L_{X}$
is also convex. Since $X_{i}:\Omega \rightarrow 2^{Y}$ is integrably bounded
and has convex weakly compact values, by Diestel's Theorem (see Appendex) it
follows that $L_{X_{i}}$ is a weakly compact subset of $L_{1}(\mu ,Y)$. More
over, $L_{X}$ is weakly compact. $L_{1}(\mu ,Y)$ equipped with the weak
topology is a locally convex topological vector space.

The correspondence $\Phi _{i}$ is convex valued, by Lemma 2 it has a
measurable graph and for each $\omega \in \Omega ,$ $\Phi _{i}(\omega ,$%
\textit{\textperiodcentered }$)$ has weakly open lower sections. Let $%
U_{i}=\{(\omega ,\widetilde{x})\in \Omega \times L_{X}:\Phi _{i}(\omega ,%
\widetilde{x})\neq \emptyset \}.$ For each $\widetilde{x}\in L_{X},$\ let%
\textit{\ }$U_{i}^{\widetilde{x}}=\{\omega \in \Omega :\Phi _{i}(\omega ,%
\widetilde{x})\neq \emptyset \}$ and for each $\omega \in \Omega ,$\ let $%
U_{i}^{\omega }=\{\widetilde{x}\in L_{X}:\Phi _{i}(\omega ,\widetilde{x}%
)\neq \emptyset \}$. The values of $\Phi _{i/U_{i}}$ have non-empty
interiors in the relative norm topology of $\ X_{i}(\omega ).$ By the
Caratheodory-type selection theorem (see Appendix), there exists a function $%
f_{i}:U_{i}\rightarrow Y$ such that $f_{i}(\omega ,\widetilde{x})\in \Phi
_{i}(\omega ,\widetilde{x})$ for all $(\omega ,\widetilde{x})\in U_{i}$ ,
for each $\widetilde{x}\in L_{X},$ $f_{i}($\textperiodcentered $,\widetilde{x%
})$ is measurable on $U_{i}^{\widetilde{x}}$, for each $\omega \in \Omega ,$ 
$f_{i}(\omega ,$\textperiodcentered $)$ is continuous on $U_{i}^{\omega }$
and moreover $f_{i}($\textperiodcentered $,$\textperiodcentered $)$ is
jointly measurable$.$

Define $G_{i}:\Omega \times L_{X}\rightarrow 2^{Y}$ by $G_{i}(\omega ,%
\widetilde{x})=\left\{ 
\begin{array}{c}
\{f_{i}(\omega ,\widetilde{x})\}\text{ if }(\omega ,\widetilde{x})\in U_{i}
\\ 
\text{cl}(A_{i}(\omega ,\widetilde{x}))_{a_{i}(\widetilde{x})}\text{ if }%
(\omega ,\widetilde{x})\notin U_{i}.%
\end{array}%
\right. $

Define $G_{i}^{^{\prime }}:L_{X}\rightarrow 2^{L_{X_{i}}}$, by $%
G_{i}^{^{\prime }}(\widetilde{x})=\{y_{i}\in L_{X_{i}}:y_{i}(\omega )\in
G_{i}(\omega ,\widetilde{x})$ $\mu -a.e.\}$ and $G^{^{\prime
}}:L_{X}\rightarrow 2^{L_{X}}$ by $G^{^{\prime }}(\widetilde{x}):=\underset{%
i\in I}{\prod }G_{i}^{^{\prime }}(\widetilde{x})$ for each $\widetilde{x}\in
L_{X}$. We shall prove that $G^{^{\prime }}$ is an upper semicontinuous
correspondence with respect to the weakly topology of $L_{X}$ and has
non-empty convex closed values. By applying Fan-Glicksberg's fixed-point
theorem [11] to $G^{^{\prime }},$ we obtain a fixed point which is the
equilibrium point for the abstract economy.

It follows by Theorem III.40 in [5] and the projection theorem\emph{\ }that
for each $\widetilde{x}\in L_{X}$, the correspondence $\widetilde{x}%
\rightarrow $cl$A_{i}(\cdot ,\widetilde{x}))_{a_{i}(\widetilde{x})}:\Omega
\rightarrow 2^{Y}$ has a measurable graph. For each $\widetilde{x}\in L_{X}$%
, the correspondence $G_{i}($\textperiodcentered $,\widetilde{x})$ has a
measurable graph. Since $\Phi _{i}(\omega ,$\textit{\textperiodcentered }$)$%
\textit{\ }has weakly open lower sections for each\textit{\ }$\omega \in
\Omega ,$ it follows that $U_{i}^{\omega }$ is weakly open in $L_{X}.$ By
Lemma 1, for each $\omega \in \Omega ,$ $G_{i}(\omega ,$\textit{%
\textperiodcentered }$):L_{X}\rightarrow 2^{Y}$ is upper semi-continuous in
the sense that the set $\{\widetilde{x}\in L_{X}:G_{i}(\omega ,\widetilde{x}%
)\}\subset V$\ is weakly open in $L_{X}$\ for every norm open subset $V$\ of 
$Y$. Moreover, $G_{i}$ is convex and non-empty valued.

$G_{i}$ is nonempty valued and for each $\widetilde{x}\in L_{X}$, $G_{i}($%
\textit{\textperiodcentered }$,\widetilde{x})$\textit{\ }has measurable
graph. Hence, by the Aumann measurable selection theorem for each fixed $%
\widetilde{x}\in L_{X}$, there exists an $\tciFourier _{i}-$measurable
function $y_{i}:\Omega \rightarrow Y$ such that $y_{i}(\omega )\in
G_{i}(\omega ,\widetilde{x})$ $\mu -a.e$. $\func{Si}$nce for each $(\omega ,%
\widetilde{x})\in \Omega \times L_{X},$ $G_{i}(\omega ,\widetilde{x})$ is
contained in the integrably bounded correspondence $X_{i}($\textit{%
\textperiodcentered }$)$, then $y_{i}\in L_{X_{i}}$ and we conclude that $%
y_{i}\in G_{i}^{^{\prime }}(\widetilde{x})$ for each $\widetilde{x}\in L_{X}$%
. Thus, $G_{i}^{^{\prime }}$ is non-empty valued.

Since for each $\widetilde{x}\in L_{X}$, $G_{i}($\textit{\textperiodcentered 
}$,\widetilde{x})$\textit{\ }has measurable graph and for each $\omega \in
\Omega ,$ $G_{i}(\omega ,$\textperiodcentered $):L_{X}\rightarrow 2^{Y}$ is
upper semicontinuous and $G_{i}(\omega ,\widetilde{x})\subset X_{i}(\omega )$
for each $(\omega ,\widetilde{x})\in \Omega \times L_{X}$, by u. s. c.
Lifting Theorem (see Appendix) it follows that $G_{i}^{^{\prime }}$ is
weakly upper semicontinuous. $G_{i}^{^{\prime }}$ is convex valued since $%
G_{i}$ is so.

$G^{^{\prime }}$ is an weakly upper semicontinuous correspondence and has
also non-empty convex closed values.

The set $L_{X}$ is weakly compact and convex, and then, by Fan-Glicksberg's
fixed-point theorem in [11], there exists $\widetilde{x}^{\ast }\in L_{X}$
such that $\widetilde{x}^{\ast }\in G^{^{\prime }}(\widetilde{x}^{\ast })$,
i.e., for each $i\in I$, $\widetilde{x}_{i}^{\ast }\in G_{i}^{^{\prime }}(%
\widetilde{x}^{\ast })$.

Then, $\widetilde{x}_{i}^{\ast }\in L_{X_{i}}$ and $\widetilde{x}_{i}^{\ast
}(\omega )\in G_{i}(\omega ,\widetilde{x}^{\ast })$ $\mu -a.e.$ Since $%
\widetilde{x}_{i}^{\ast }(\omega )\notin (A_{i}(\omega ,\widetilde{x}^{\ast
}))_{a_{i}(\widetilde{x}^{\ast })}\cap (P_{i}(\omega ,\widetilde{x}^{\ast
}))_{p_{i}(\widetilde{x}^{\ast })}$ $\mu -a.e,$ it follows that $(\omega ,%
\widetilde{x}^{\ast })\notin U_{i}$ for each $i\in I$ and $\widetilde{x}%
_{i}^{\ast }\in $cl$(A_{i}(\omega ,\widetilde{x}^{\ast }))_{a_{i}(\widetilde{%
x}^{\ast })}$ $\mu -a.e$. We have also that $(A_{i}(\omega ,\widetilde{x}%
^{\ast }))_{a_{i}(\widetilde{x}^{\ast })}\cap (P_{i}(\omega ,\widetilde{x}%
^{\ast }))_{p_{i}(\widetilde{x}^{\ast })}=\emptyset .\Box \medskip $

\begin{theorem}
\textit{Let }$I$\textit{\ be a countable or uncounatble set. Let the family} 
$G=\{(\Omega ,\tciFourier ,\mu ),$ $(X_{i},\tciFourier
_{i},A_{i},P_{i},a_{i},b_{i},z_{i})_{i\in I},\}$ \textit{be a general
Bayesian abstract economy satisfying A.1)-A.4). Then there exists a Bayesian
fuzzy equilibrium for }$G$\textit{.}
\end{theorem}

\textit{For each }$i\in I:$

\textit{A.1)}

\textit{\ \ \ \ \ (a) }$X_{i}:\Omega \rightarrow \mathcal{F}(Y)$\textit{\ is
such that }$\omega \rightarrow (X_{i}(\omega ))_{z_{i}}:\Omega \rightarrow
2^{Y}$ \textit{is} \textit{a nonempty, convex, weakly compact-valued and
integrably bounded correspondence;}

\textit{\ \ \ \ \ (b) }$X_{i}:\Omega \rightarrow \mathcal{F}(Y)$\textit{\ is
such that }$\omega \rightarrow (X_{i}(\omega ))_{z_{i}}:\Omega \rightarrow
2^{Y}$ \textit{is }$\tciFourier _{i}-$\textit{lower measurable;}

\textit{\ A.2)}

\textit{\ \ \ \ \ (a) For each }$(\omega ,\widetilde{x})\in \Omega \times
L_{X},$\textit{\ }$(A_{i}(\omega ,\widetilde{x}))_{a_{i}(\widetilde{x})}$%
\textit{\ is\ nonempty convex and compact}$.$\textit{\medskip }

\textit{\ \ \ \ \ (b)\ For each }$\widetilde{x}\in L_{X},$ \textit{the
correspondence }$\omega \rightarrow (A_{i}(\omega ,\widetilde{x}))_{a_{i}(%
\widetilde{x})}:\Omega \rightarrow 2^{Y}$\textit{\ has measurable graph;}

\textit{\ \ \ \ (c) For each }$\omega \in \Omega ,$\textit{\ }$\widetilde{x}%
\rightarrow (A_{i}(\omega ,\widetilde{x}))_{a_{i}(\widetilde{x}%
)}:L_{X}\rightarrow 2^{Y}$\textit{\ is upper semicontinuous in the sense
that the set }$\{\widetilde{x}\in L_{X}\mathit{:}(A_{i}(\omega ,\widetilde{x}%
))_{a_{i}(\widetilde{x})}\}\subset V$ \textit{is weakly open in }$L_{X}$%
\textit{\ for every norm open subset }$V$\textit{\ of} $Y$;

\textit{A.3)}

\textit{\ \ \ \ \ (a) For each }$(\omega ,\widetilde{x})\in \Omega \times
L_{X},$ \textit{the correspondence }$\omega \rightarrow (P_{i}(\omega ,%
\widetilde{x}))_{p_{i}(\widetilde{x})}:\Omega \times L_{X}\rightarrow 2^{Y}$%
\textit{\ has nonempty convex compact values such that }$(P_{i}(\omega ,%
\widetilde{x}))_{p_{i}(\widetilde{x})}\subset (X(\omega ))_{z_{i}}$\textit{\
for each }$(\omega ,\widetilde{x})\in \Omega \times L_{X};$

\textit{\ \ \ \ (b)\ For each }$\widetilde{x}\in L_{X},$ \textit{the
correspondence} $\omega \rightarrow (P_{i}(\omega ,\widetilde{x}))_{p_{i}(%
\widetilde{x})}:\Omega \times L_{X}\rightarrow 2^{Y}$ \textit{has a
measurable graph; \ }

\textit{(c) For each }$\omega \in \Omega ,$\textit{\ }$\widetilde{x}%
\rightarrow (P_{i}(\omega ,\widetilde{x}))_{p_{i}(\widetilde{x}%
)}:L_{X}\rightarrow 2^{Y}$\textit{\ is upper semicontinuous}$;$

\textit{A.4)}

\textit{\ \ \ \ (a) For each }$\widetilde{x}\in L_{X},$\textit{\ for each }$%
\omega \in \Omega ,$\textit{\ }$\widetilde{x}_{i}(\omega )\notin
(A_{i}(\omega ,\widetilde{x}))_{a_{i}(\widetilde{x})}\cap (P_{i}(\omega ,%
\widetilde{x}))_{p_{i}(\widetilde{x})};$

\ \ \ (b) \textit{For each }$\omega \in \Omega ,$ \textit{the set }$%
U_{i}^{\omega }=\{\widetilde{x}\in L_{X}:$\textit{\ }$(A_{i}(\omega ,%
\widetilde{x}))_{a_{i}(\widetilde{x})}\cap (P_{i}(\omega ,\widetilde{x}%
))_{p_{i}(\widetilde{x})}\neq \emptyset \}$\textit{\ is weakly open in }$%
L_{X}.\medskip $

\textit{Proof.} For each $i\in I,$ define $\Phi _{i}:\Omega \times
L_{X}\rightarrow 2^{Y}$ by $\Phi _{i}(\omega ,\widetilde{x})=(A_{i}(\omega ,%
\widetilde{x}))_{a_{i}(\widetilde{x})}\cap (P_{i}(\omega ,\widetilde{x}%
))_{p_{i}(\widetilde{x})}.$

We prove first that $L_{X}$ is a non-empty, convex, weakly compact subset in 
$L_{1}(\mu ,Y).$

Since $(\Omega $\textit{, }$\tciFourier ,\mu )$ is a complete finite measure
space, $Y$ is a separable Banach space and $(X_{i}(\cdot ))_{z_{i}}:\Omega
\rightarrow 2^{Y}$ has measurable graph, by Aumann's selection theorem (see
Appendex) it follows that there exists a $\tciFourier _{i}$-measurable
function $f_{i}:\Omega \rightarrow Y$ such that $f_{i}(\omega )\in
X_{i}(\omega )$ $\mu -a.e$. Since $(X_{i}(\cdot ))_{z_{i}}$ is integrably
bounded, we have that $f_{i}\in L_{1}(\mu ,Y)$, hence $L_{X_{i}}$ is
non-empty and $L_{X}=\tprod\limits_{i\in I}L_{X_{i}}$ is non-empty.
Obviously $L_{X_{i}}$ is convex and $L_{X}$ is also convex. Since $%
X_{i}:\Omega \rightarrow 2^{Y}$ is integrably bounded and has convex weakly
compact values, by Diestel's Theorem (see Appendex) it follows that $%
L_{X_{i}}$ is a weakly compact subset of $L_{1}(\mu ,Y)$. More over, $L_{X}$
is weakly compact. $L_{1}(\mu ,Y)$ equipped with the weak topology is a
locally convex topological vector space.

Define $G_{i}:\Omega \times L_{X}\rightarrow 2^{Y}$ by $G_{i}(\omega ,%
\widetilde{x})=\left\{ 
\begin{array}{c}
\Phi _{i}(\omega ,\widetilde{x})\text{ \ \ \ if \ \ \ }(\omega ,\widetilde{x}%
)\in U_{i}; \\ 
(A_{i}(\omega ,\widetilde{x}))_{a_{i}(\widetilde{x})}\text{ if }(\omega ,%
\widetilde{x})\notin U_{i}.%
\end{array}%
\right. $

For each $\widetilde{x}\in L_{X},$ the correspondence $\omega \rightarrow
(A_{i}(\omega ,\widetilde{x}))_{a_{i}(\widetilde{x})}:\Omega \rightarrow
2^{Y}$ has a measurable graph. Hence for each $\widetilde{x}\in L_{X}$, the
correspondence $G_{i}($\textperiodcentered $,\widetilde{x})$ has a
measurable graph. By the assumption A4) (b), we have that $U_{i}^{\omega }$
is weakly open in $L_{X}.$ By Lemma 1, for each $\omega \in \Omega ,$ $%
G_{i}(\omega ,$\textit{\textperiodcentered }$):L_{X}\rightarrow 2^{Y}$ is
upper semi-continuous. Moreover, $G_{i}$ is non-empty convex compact valued.

Define $G_{i}^{^{\prime }}:L_{X}\rightarrow 2^{L_{X_{i}}}$, by $%
G_{i}^{^{\prime }}(\widetilde{x})=\{y_{i}\in L_{X_{i}}:y_{i}(\omega )\in
G_{i}(\omega ,\widetilde{x})$ $\mu -a.e.\}$ and $G^{^{\prime
}}:L_{X}\rightarrow 2^{L_{X}}$ by $G^{^{\prime }}(\widetilde{x}):=\underset{%
i\in I}{\prod }G_{i}^{^{\prime }}(\widetilde{x})$ for each $x\in L_{X}$.

Since for each $\widetilde{x}\in L_{X}$, $G_{i}($\textit{\textperiodcentered 
}$,\widetilde{x})$\textit{\ }has a measurable graph and for each $\omega \in
\Omega ,$ $G_{i}(\omega ,$\textperiodcentered $):L_{X}\rightarrow 2^{Y}$ is
upper semicontinuous and $G_{i}(\omega ,\widetilde{x})\subset (X_{i}(\omega
))_{z_{i}}$ for each $(\omega ,\widetilde{x})\in \Omega \times L_{X}$, by u.
s. c. Lifting Theorem (see Appendix) it follows that $G_{i}^{^{\prime }}$ is
weakly upper semicontinuous. $G_{i}^{^{\prime }}$ is convex valued since $%
G_{i}$ is so.

$G^{^{\prime }}$ is an weakly upper semicontinuous correspondence and has
also non-empty convex closed values.

The set $L_{X}$ is weakly compact and convex, and then, by Fan-Glicksberg's
fixed-point theorem (1952), there exists $\widetilde{x}^{\ast }\in L_{X}$
such that $\widetilde{x}^{\ast }\in G^{^{\prime }}(x^{\ast })$, i.e., for
each $i\in I$, $\widetilde{x}_{i}^{\ast }\in G_{i}^{^{\prime }}(\widetilde{x}%
^{\ast })$.

Then, $\widetilde{x}_{i}^{\ast }\in L_{X_{i}}$ and $\widetilde{x}_{i}^{\ast
}(\omega )\in G_{i}(\omega ,\widetilde{x}^{\ast })$ $\mu -a.e.$ Since $%
\widetilde{x}_{i}^{\ast }(\omega )\notin A_{i}(\omega ,\widetilde{x}^{\ast
})\cap P_{i}(\omega ,\widetilde{x}^{\ast })$ $\mu -a.e,$ it follows that $%
(\omega ,\widetilde{x}^{\ast })\notin U_{i}$ for each $i\in I$ and $%
\widetilde{x}_{i}^{\ast }\in $cl$A_{i}(\widetilde{x}^{\ast })$ $\mu -a.e$.
We have also that $A_{i}(\omega ,\widetilde{x}^{\ast })\cap P_{i}(\omega ,%
\widetilde{x}^{\ast })=\emptyset .\Box \medskip $

\section{RANDOM\ QUASI-VARIATIONAL INEQUALITIES}

In this section, we are establishing new random quasi variational
inequalities with random fuzzy mappings and random fixed point theorems. The
proofs rely on the theorem of Bayesian fuzzy equilibrium existence for the
Bayesian abstract fuzzy economy.

This is our first theorem.

\begin{theorem}
\textit{Let }$I$\textit{\ be a countable or uncounatble set. Let }$(\Omega $%
\textit{, }$\tciFourier ,\mu )$ be a complete finite separable measure space
and let $Y$ be a separable Banach space.\textit{\ Suppose that the following
conditions are satisfied:}
\end{theorem}

\textit{For each }$i\in I:$

\textit{A.1)}

\textit{\ \ \ \ \ (a) }$X_{i}:\Omega \rightarrow \mathcal{F}(Y)$\textit{\ is
such that }$\omega \rightarrow X_{i}(\omega )_{z}:\Omega \rightarrow 2^{Y}$ 
\textit{is} \textit{a nonempty, convex, weakly compact-valued and integrably
bounded correspondence.}

\textit{\ \ \ \ \ (b) }$X_{i}:\Omega \rightarrow \mathcal{F}(Y)$\textit{\ is
such that }$\omega \rightarrow (X_{i}(\omega ))_{z}:\Omega \rightarrow 2^{Y}$
\textit{is }$\tciFourier _{i}-$\textit{lower measurable;}

\textit{\ A.2)}

\textit{\ \ \ \ \ (a) For each }$(\omega ,\widetilde{x})\in \Omega \times
L_{X},$\textit{\ }$(A_{i}(\omega ,\widetilde{x}))_{a_{i}(\widetilde{x})}$%
\textit{\ is\ convex and has a non-empty interior in the relative norm
topology of }$(X_{i}(\omega ))_{z}.$\textit{\medskip }

\textit{\ \ \ \ \ (b)\ the correspondence }$(\omega ,\widetilde{x}%
)\rightarrow (A_{i}(\omega ,\widetilde{x}))_{a_{i}(\widetilde{x})}:\Omega
\times L_{X}\rightarrow 2^{Y}$\textit{\ has measurable graph i.e. }$%
\{(\omega ,\widetilde{x},y)\in \Omega \times L_{X}\times Y:y\in
(A_{i}(\omega ,\widetilde{x}))_{a_{i}(\widetilde{x})}\}\in \tciFourier
\otimes $\textit{\ss }$_{w}(L_{X})\otimes $\textit{\ss }$(Y)$\textit{\ where 
\ss }$_{w}(L_{X})$\textit{\ is the Borel }$\sigma -$\textit{algebra for the
weak topology on }$L_{X}$\textit{\ and \ss }$(Y)$\textit{\ is the Borel }$%
\sigma -$\textit{algebra for the norm topology on }$Y$\textit{.}

\textit{\ \ \ \ \ (c)\ the correspondence }$(\omega ,\widetilde{x}%
)\rightarrow (A_{i}(\omega ,\widetilde{x}))_{a_{i}(\widetilde{x})}$\textit{\
has weakly open lower sections, i.e., for each }$\omega \in \Omega $\textit{%
\ and for each }$y\in Y,$\textit{\ the set }$((A_{i}(\omega ,\widetilde{x}%
)_{a_{i}(\widetilde{x})})^{-1}(\omega ,y)=\{\widetilde{x}\in L_{X}:y\in
(A_{i}(\omega ,\widetilde{x}))_{a_{i}(\widetilde{x})}\}\}$\textit{\ is
weakly open in }$L_{X};$

\textit{\ \ \ \ \ (d) For each }$\omega \in \Omega ,$\textit{\ }$\widetilde{x%
}\rightarrow $cl$(A_{i}(\omega ,\widetilde{x}))_{a_{i}(\widetilde{x}%
)}:L_{X}\rightarrow 2^{Y}$\textit{\ is upper semicontinuous in the sense
that the set }$\{\widetilde{x}\in L_{X}\mathit{:}$cl$(A_{i}(\omega ,%
\widetilde{x}))_{a_{i}(\widetilde{x})}\subset V\}$ \textit{is weakly open in 
}$L_{X}$\textit{\ for every norm open subset }$V$\textit{\ of }$Y$\textit{.}

\textit{A.3)}

\ \ \ \textit{\ }$\psi _{i}:\Omega \times L_{X}\times Y\rightarrow R\cup
\{-\infty ,+\infty \}$\textit{\ is such that:}

\ \ \ (\textit{a) }$\widetilde{x}\rightarrow \psi _{i}(\omega ,\widetilde{x}%
,y)$\textit{\ is lower semicontinuous on }$L_{X}$\textit{\ for each fixed }$%
(\omega ,y)\in \Omega \times Y;$

\ \ \ (\textit{b) }$\widetilde{x}(\omega )\notin $cl$\{y\in Y:\psi
_{i}(\omega ,\widetilde{x},y)>0\}$\textit{\ for each fixed }$(\omega ,%
\widetilde{x})\in \Omega \times L_{X};$

\ \ \ (\textit{c) for each} $(\omega ,\widetilde{x})\in \Omega \times L_{X},$
$\psi _{i}(\omega ,\widetilde{x},\cdot )$ \textit{is concave;}

\ \ \ (\textit{d) for each }$\omega \in \Omega ,$\textit{\ }$\{\widetilde{x}%
\in L_{X}:\alpha _{i}(\omega ,\widetilde{x})>0\}$\textit{\ is weakly open in}
$L_{X},$\textit{\ where }$\alpha _{i}:\Omega \times L_{X}\rightarrow R$%
\textit{\ is defined by }$\alpha _{i}(\omega ,\widetilde{x})=\sup_{y\in
(A_{i}(\omega ,\widetilde{x}))_{a_{i}(\widetilde{x})}}\psi _{i}(\omega ,%
\widetilde{x},y)$\textit{\ for each }$(\omega ,\widetilde{x})\in \Omega
\times L_{X};$

\ \ \ (\textit{e) }$\{(\omega ,\widetilde{x}):\alpha _{i}(\omega ,\widetilde{%
x})>0\}\in \mathcal{F}_{i}\otimes B(L_{X})$\textit{.}

\textit{Then, there exists }$\widetilde{x}^{\ast }\in L_{X}$\textit{\ such
that }for every $i\in I$,

i) $\widetilde{x_{i}}^{\ast }(\omega )\in $cl$(A_{i}(\omega ,\widetilde{x}%
^{\ast }))_{a_{i}(\widetilde{x}^{\ast })};$

\textit{ii) sup}$_{y\in (A_{i}(\omega ,\widetilde{x}^{\ast }))_{a_{i}(%
\widetilde{x}^{\ast })}}\psi _{i}(\omega ,\widetilde{x}^{\ast },y)\leq 0$%
\textit{.}$\medskip $

\textit{Proof.} For every $i\in I,$ let $P_{i}:\Omega \times
S_{X}^{1}\rightarrow \mathcal{F}(Y)$ and $p_{i}:S_{X}^{1}\rightarrow (0,1]$
such that $(P_{i}(\omega ,\widetilde{x}))_{p_{i}(\widetilde{x})}=\{y\in
Y:\psi _{i}(\omega ,\widetilde{x},y)>0\}$ for each $(\omega ,\widetilde{x}%
)\in \Omega \times L_{X}.$

We shall show that the abstract economy $G=\{(\Omega ,\tciFourier ,\mu ),$ $%
(X_{i},\tciFourier _{i},A_{i},P_{i}a_{i},p_{i},z_{z})_{i\in I},\}$ satisfies
all hypotheses of Theorem 1.

Suppose $\omega \in \Omega .$

According to A3 a), we have that\textit{\ }$\widetilde{x}\rightarrow
(P_{i}(\omega ,\widetilde{x}))_{p_{i}(\widetilde{x})}:\Omega \rightarrow
2^{Y}$\textit{\ }has open lower sections with nonempty compact values and
according to A3 b), $\widetilde{x_{i}}(\omega )\not\in (P_{i}(\omega ,%
\widetilde{x}))_{p_{i}(\widetilde{x})}$ for each $\widetilde{x}\in L_{X}.$
Assumption A3 c) implies that $\widetilde{x}\rightarrow (P_{i}(\omega ,%
\widetilde{x}))_{p_{i}(\widetilde{x})}:\Omega \rightarrow 2^{Y}$ has convex
values.

By the definition of $\alpha _{i},$ we note that $\{\widetilde{x}\in
L_{X}:(A_{i}(\omega ,\widetilde{x}))_{a_{i}(\widetilde{x})}\cap
(P_{i}(\omega ,\widetilde{x}))_{p_{i}(\widetilde{x})}\neq \emptyset \}==\{%
\widetilde{x}\in L_{X}:\alpha _{i}(\omega ,\widetilde{x})>0\}$ so that $\{%
\widetilde{x}\in L_{X}:(A_{i}(\omega ,\widetilde{x}))_{a_{i}(\widetilde{x}%
)}\cap (P_{i}(\omega ,\widetilde{x}))_{p_{i}(\widetilde{x})}\neq \emptyset
\} $ is weakly open in $L_{X}$ by A3 c).

According to A2 b) and A3 e), it follows that the correspondences\textit{\ }$%
(\omega ,\widetilde{x})\rightarrow (A_{i}(\omega ,\widetilde{x}))_{a_{i}(%
\widetilde{x})}:\Omega \times L_{X}\rightarrow 2^{Y}$ and $(\omega ,%
\widetilde{x})\rightarrow (P_{i}(\omega ,\widetilde{x}))_{p_{i}(\widetilde{x}%
)}:\Omega \times L_{X}\rightarrow 2^{Y}$\textit{\ }have measurable graphs$.$

Thus the Bayesian abstract fuzzy economy $G=\{(\Omega ,\tciFourier ,\mu ),$ $%
(X_{i},\tciFourier _{i},A_{i},P_{i},a_{i},b_{i},z_{i})_{i\in I},\}$
satisfies all hypotheses of Theorem 1. Therefore, there exists $\widetilde{x}%
^{\ast }\in L_{X}$ such that for every $i\in I:$

$\widetilde{x}_{i}^{\ast }(\omega )\in $cl$(A_{i}(\omega ,\widetilde{x}%
^{\ast }))_{a_{i}(\widetilde{x}^{\ast })}$ $\mu -a.e$ and

$(A_{i}(\omega ,\widetilde{x}^{\ast }))_{a_{i}(\widetilde{x}^{\ast })}\cap
(P_{i}(\omega ,\widetilde{x}^{\ast }))_{p_{i}(\widetilde{x}^{\ast })}=\phi $ 
$\mu -a.e;$

that is, there exists $\widetilde{x}^{\ast }\in L_{X}$ such that for every $%
i\in I:$

i) $\widetilde{x_{i}}^{\ast }(\omega )\in $cl$(A_{i}(\omega ,\widetilde{x}%
^{\ast }))_{a_{i}(\widetilde{x}^{\ast })};$

ii) sup$_{y\in (A_{i}(\omega ,\widetilde{x}^{\ast }))_{a_{i}(\widetilde{x}%
^{\ast })}}\psi _{i}(\omega ,\widetilde{x}^{\ast },y)\leq 0$.$\medskip $

If \TEXTsymbol{\vert}I\TEXTsymbol{\vert}=1, we obtain the following
corollary.

\begin{corollary}
\textit{Let }$(\Omega $\textit{, }$\tciFourier ,\mu )$ be a complete finite
separable measure space and let $Y$ be a separable Banach space.\textit{\
Suppose that the following conditions are satisfied:}
\end{corollary}

\textit{A.1)}

\textit{\ \ \ \ \ (a) }$X:\Omega \rightarrow \mathcal{F}(Y)$\textit{\ is
such that }$\omega \rightarrow (X(\omega ))_{z}:\Omega \rightarrow 2^{Y}$ 
\textit{is} \textit{a nonempty, convex, weakly compact-valued and integrably
bounded correspondence.}

\textit{\ \ \ \ \ (b) }$X:\Omega \rightarrow \mathcal{F}(Y)$\textit{\ is
such that }$\omega \rightarrow (X(\omega ))_{z}:\Omega \rightarrow 2^{Y}$ 
\textit{is }$\tciFourier -$\textit{lower measurable;}

\textit{\ A.2)}

\textit{\ \ \ \ \ (a) For each }$(\omega ,\widetilde{x})\in \Omega \times
L_{X},$\textit{\ }$(A(\omega ,\widetilde{x}))_{a(\widetilde{x})}$\textit{\
is\ convex and has a non-empty interior in the relative norm topology of }$%
(X(\omega ))_{z}.$\textit{\medskip }

\textit{\ \ \ \ \ (b)\ the correspondence }$(\omega ,\widetilde{x}%
)\rightarrow (A(\omega ,\widetilde{x}))_{a_{i}(\widetilde{x})}:\Omega \times
L_{X}\rightarrow 2^{Y}$\textit{\ has measurable graph i.e. }$\{(\omega ,%
\widetilde{x},y)\in \Omega \times L_{X}\times Y:y\in (A(\omega ,\widetilde{x}%
))_{a(\widetilde{x})}\}\in \tciFourier \otimes $\textit{\ss }$%
_{w}(L_{X})\otimes $\textit{\ss }$(Y)$\textit{\ where \ss }$_{w}(L_{X})$%
\textit{\ is the Borel }$\sigma -$\textit{algebra for the weak topology on }$%
L_{X}$\textit{\ and \ss }$(Y)$\textit{\ is the Borel }$\sigma -$\textit{%
algebra for the norm topology on }$Y$\textit{.}

\textit{\ \ \ \ \ (c)\ the correspondence }$(\omega ,\widetilde{x}%
)\rightarrow (A(\omega ,\widetilde{x}))_{a(\widetilde{x})}$\textit{\ has
weakly open lower sections, i.e., for each }$\omega \in \Omega $\textit{\
and for each }$y\in Y,$\textit{\ the set }$((A(\omega ,\widetilde{x})_{a(%
\widetilde{x})})^{-1}(\omega ,y)=\{\widetilde{x}\in L_{X}:y\in (A(\omega ,%
\widetilde{x}))_{a(\widetilde{x})}\}\}$\textit{\ is weakly open in }$L_{X};$

\textit{\ \ \ \ \ (d) For each }$\omega \in \Omega ,$\textit{\ }$\widetilde{x%
}\rightarrow $cl$(A(\omega ,\widetilde{x}))_{a(\widetilde{x}%
)}:L_{X}\rightarrow 2^{Y}$\textit{\ is upper semicontinuous in the sense
that the set }$\{\widetilde{x}\in L_{X}\mathit{:}$cl$(A(\omega ,\widetilde{x}%
))_{a(\widetilde{x})}\subset V\}$ \textit{is weakly open in }$L_{X}$\textit{%
\ for every norm open subset }$V$\textit{\ of }$Y$\textit{.}

\textit{A.3)}

\ \ \ \textit{\ }$\psi _{i}:\Omega \times L_{X}\times Y\rightarrow R\cup
\{-\infty ,+\infty \}$\textit{\ is such that:}

\ \ \ (\textit{a) }$\widetilde{x}\rightarrow \psi (\omega ,\widetilde{x},y)$%
\textit{\ is lower semicontinuous on }$L_{X}$\textit{\ for each fixed }$%
(\omega ,y)\in \Omega \times Y;$

\ \ \ (\textit{b) }$\widetilde{x}(\omega )\notin $cl$\{y\in Y:\psi (\omega ,%
\widetilde{x},y)>0\}$\textit{\ for each fixed }$(\omega ,\widetilde{x})\in
\Omega \times L_{X};$

\ \ \ (\textit{c) for each} $(\omega ,\widetilde{x})\in \Omega \times L_{X},$
$\psi (\omega ,\widetilde{x},\cdot )$ \textit{is concave;}

\ \ \ (\textit{d) for each }$\omega \in \Omega ,$\textit{\ }$\{\widetilde{x}%
\in L_{X}:\alpha (\omega ,\widetilde{x})>0\}$\textit{\ is weakly open in} $%
L_{X},$\textit{\ where }$\alpha :\Omega \times L_{X}\rightarrow R$\textit{\
is defined by }$\alpha (\omega ,\widetilde{x})=\sup_{y\in (A(\omega ,%
\widetilde{x}))_{a(\widetilde{x})}}\psi (\omega ,\widetilde{x},y)$\textit{\
for each }$(\omega ,\widetilde{x})\in \Omega \times L_{X};$

\ \ \ (\textit{e) }$\{(\omega ,\widetilde{x}):\alpha (\omega ,\widetilde{x}%
)>0\}\in \mathcal{F}\otimes B(L_{X})$\textit{.}

\textit{Then, there exists }$\widetilde{x}^{\ast }\in L_{X}$\textit{\ such
that}:

i) $\widetilde{x}^{\ast }(\omega )\in $cl$(A(\omega ,\widetilde{x}^{\ast
}))_{a(\widetilde{x}^{\ast })};$

\textit{ii) sup}$_{y\in (A(\omega ,\widetilde{x}^{\ast }))_{a(\widetilde{x}%
^{\ast })}}\psi (\omega ,\widetilde{x}^{\ast },y)\leq 0$\textit{.}$\medskip $

As a consequence of Theorem 2, we prove the following random generalized
quasi-variational inequality with random fuzzy mappings. This theorem is
comparable with Theorem 4.1 in [27], which is valid in a non-fuzzy framework
and concerns upper-semicontinuos correspondences defined on metrizable
spaces.

\begin{theorem}
\textit{Let }$(\Omega $\textit{, }$\tciFourier ,\mu )$ be a complete finite
separable measure space and let $Y$ be a separable Banach space.\textit{\
Suppose that the following conditions are satisfied:}
\end{theorem}

\textit{For each }$i\in I:$

\textit{A.1)}

\textit{\ \ \ \ \ (a) }$X_{i}:\Omega \rightarrow \mathcal{F}(Y)$\textit{\ is
such that }$\omega \rightarrow (X_{i}(\omega ))_{z}:\Omega \rightarrow 2^{Y}$
\textit{is} \textit{a nonempty, convex, weakly compact-valued and integrably
bounded correspondence.}

\textit{\ \ \ \ \ (b) }$X_{i}:\Omega \rightarrow \mathcal{F}(Y)$\textit{\ is
such that }$\omega \rightarrow (X_{i}(\omega ))_{z}:\Omega \rightarrow 2^{Y}$
\textit{is }$\tciFourier _{i}-$\textit{lower measurable;}

\textit{\ A.2)}

\textit{\ \ \ \ \ (a) For each }$(\omega ,\widetilde{x})\in \Omega \times
L_{X},$\textit{\ }$(A_{i}(\omega ,\widetilde{x}))_{a_{i}(\widetilde{x})}$%
\textit{\ is\ convex and has a non-empty interior in the relative norm
topology of }$(X_{i}(\omega ))_{z}.$\textit{\medskip }

\textit{\ \ \ \ \ (b)\ the correspondence }$(\omega ,\widetilde{x}%
)\rightarrow (A_{i}(\omega ,\widetilde{x}))_{a_{i}(\widetilde{x})}:\Omega
\times L_{X}\rightarrow 2^{Y}$\textit{\ has measurable graph i.e. }$%
\{(\omega ,\widetilde{x},y)\in \Omega \times L_{X}\times Y:y\in
(A_{i}(\omega ,\widetilde{x}))_{a_{i}(\widetilde{x})}\}\in \tciFourier
\otimes $\textit{\ss }$_{w}(L_{X})\otimes $\textit{\ss }$(Y)$\textit{\ where 
\ss }$_{w}(L_{X})$\textit{\ is the Borel }$\sigma -$\textit{algebra for the
weak topology on }$L_{X}$\textit{\ and \ss }$(Y)$\textit{\ is the Borel }$%
\sigma -$\textit{algebra for the norm topology on }$Y$\textit{.}

\textit{\ \ \ \ \ (c)\ the correspondence }$(\omega ,\widetilde{x}%
)\rightarrow (A_{i}(\omega ,\widetilde{x}))_{a_{i}(\widetilde{x})}$\textit{\
has weakly open lower sections, i.e., for each }$\omega \in \Omega $\textit{%
\ and for each }$y\in Y,$\textit{\ the set }$((A_{i}(\omega ,\widetilde{x}%
)_{a_{i}(\widetilde{x})})^{-1}(\omega ,y)=\{\widetilde{x}\in L_{X}:y\in
(A_{i}(\omega ,\widetilde{x}))_{a_{i}(\widetilde{x})}\}\}$\textit{\ is
weakly open in }$L_{X};$

\textit{\ \ \ \ \ (d) For each }$\omega \in \Omega ,$\textit{\ }$\widetilde{x%
}\rightarrow $cl$(A_{i}(\omega ,\widetilde{x}))_{a_{i}(\widetilde{x}%
)}:L_{X}\rightarrow 2^{Y}$\textit{\ is upper semicontinuous in the sense
that the set }$\{\widetilde{x}\in L_{X}\mathit{:}$cl$(A_{i}(\omega ,%
\widetilde{x}))_{a_{i}(\widetilde{x})}\subset V\}$ \textit{is weakly open in 
}$L_{X}$\textit{\ for every norm open subset }$V$\textit{\ of }$Y$\textit{;}

\textit{A.3)}

$G_{i}:\Omega \times Y\rightarrow \mathcal{F}(Y^{\prime })$ \textit{and} $%
g_{i}:Y\rightarrow (0,1]$ \textit{are such that:}

\ \ \ (\textit{a) for each }$\omega \in \Omega ,$\textit{\ }$y\rightarrow
(G_{i}(\omega ,y))_{g_{i}(y)}:Y\rightarrow 2^{Y^{\prime }}$ \textit{is
monotone\ with non-empty values}$;$

\ \ \ (\textit{b) for each }$\omega \in \Omega ,$\textit{\ }$y\rightarrow
(G_{i}(\omega ,y))_{g_{i}(y)}:L\cap Y\rightarrow 2^{Y^{\prime }}$\textit{\
is lower semicontinuous from the relative topology of }$Y$\textit{\ into the
weak}$^{\ast }-$\textit{topology }$\sigma (Y^{\prime },Y)$ \textit{of }$%
Y^{\prime }$\textit{\ for each one-dimensional flat }$L\subset Y;$

\ \ \ \textit{A.4)}

\ \ \ (\textit{a) }$f_{i}:\Omega \times L_{X}\times Y\rightarrow R\cup
\{\infty ,-\infty \}$\textit{\ is such that }$\widetilde{x}\rightarrow
f_{i}(\omega ,\widetilde{x},y)$\textit{\ is lower semicontinuous on }$L_{X}$%
\textit{\ for each fixed }$(\omega ,y)\in \Omega \times Y,$\textit{\ }$%
f_{i}(\omega ,\widetilde{x},\widetilde{x}(t,\omega ))=0$\textit{\ for each }$%
(\omega ,\widetilde{x})\in \Omega \times L_{X}$\textit{\ and }$y\rightarrow
f_{i}(\omega ,\widetilde{x},y)$\textit{\ is concave on }$Y$ \textit{for each
fixed }$(\omega ,\widetilde{x})\in \Omega \times L_{X}$\textit{;}

\ \ \ (\textit{b) for each fixed }$\omega \in \Omega ,$\textit{\ the set}

\textit{\ }$\{\widetilde{x}\in S_{X}^{1}:\sup_{y\in (A_{i}(\omega ,%
\widetilde{x}))_{a_{i}(\widetilde{x})})}[\sup_{u\in (G_{i}(\omega
,y))_{g_{i}(y)}}$\textit{Re}$\langle u,\widetilde{x}-y\rangle +f_{i}(\omega ,%
\widetilde{x},y)]>0\}$\textit{\ is weakly open in }$L_{X}$

\ \ \ (c) $\{(\omega ,\widetilde{x}):\sup_{u\in (G_{i}(\omega
,y))_{g_{i}(y)}}\mathit{Re}\langle u,\widetilde{x}-y\rangle +f_{i}(\omega ,%
\widetilde{x},y)>0\}\in \mathcal{F}\otimes B(L_{X})$\textit{.}

\textit{Then, there exists }$\widetilde{x}^{\ast }\in L_{X}$\textit{\ such
that }for every $i\in I$:

i) $\widetilde{x_{i}}^{\ast }(\omega )\in $cl$(A_{i}(\omega ,\widetilde{x}%
^{\ast }))_{a_{i}(\widetilde{x}^{\ast })};$

ii) \textit{sup}$_{u\in (G_{i}(\omega ,\widetilde{x}^{\ast }(\omega
)))_{g_{i}(\widetilde{x}^{\ast }(\omega ))}}Re\langle u,\widetilde{x}^{\ast
}(\omega )-y\rangle +f_{i}(\omega ,\widetilde{x},y)]\leq 0$\textit{\ for all 
}$y\in (A_{i}(\omega ,\widetilde{x}^{\ast }))_{a_{i}(\widetilde{x}^{\ast
})}))$\textit{.}$\medskip $

\textit{Proof.} Let us define $\psi _{i}:\Omega \times L_{X}\times
Y\rightarrow R\cup \{-\infty ,+\infty \}$ by

$\psi _{i}(\omega ,\widetilde{x},y)=\sup_{u\in (G_{i}(\omega
,y))_{g_{i}(y)}} $\textit{Re}$\langle u,\widetilde{x}-y\rangle
+f_{i}(t,\omega ,\widetilde{x},y)$ for each $(\omega ,\widetilde{x},y)\in
\Omega \times L_{X}\times Y.$

According to assumption A4 a), $\widetilde{x}\rightarrow f_{i}(t,\omega ,%
\widetilde{x},y)$ is lower semicontinuous on $L_{X}$ for each fixed $(\omega
,y)\in \Omega \times Y$ and $f_{i}(\omega ,\widetilde{x},\widetilde{x}%
(t,\omega ))=0$ for each $(\omega ,\widetilde{x})\in \Omega \times L_{X}$
implies that $\widetilde{x}(\omega )\notin \{y\in Y:\psi _{i}(\omega ,%
\widetilde{x},y)>0\}$\textit{\ }for each fixed\textit{\ }$(\omega ,%
\widetilde{x})\in \Omega \times L_{X}.$

We also have that for each $(\omega ,\widetilde{x})\in \Omega \times L_{X},$ 
$\psi _{i}(\omega ,\widetilde{x},\cdot )$ is concave$.$ This fact is a
consequence of assumption A3 b).

All hypotheses of Theorem 2 are satisfied. According to Theorem 2, there
exists $\widetilde{x}^{\ast }\in L_{X}$ such that $\widetilde{x_{i}}^{\ast
}(\omega )\in $cl$(A_{i}(\omega ,\widetilde{x}^{\ast }))_{a_{i}(\widetilde{x}%
^{\ast })}$ for every $i\in I$

and

(1) \ \ sup$_{y\in A_{i}(\omega ,\widetilde{x}^{\ast }))_{a_{i}(\widetilde{x}%
^{\ast })}}\sup_{u\in (G_{i}(\omega ,y))_{g_{i}(y)}}[\mathit{Re}\langle u,%
\widetilde{x}^{\ast }(\omega )-y\rangle +f_{i}(\omega ,\widetilde{x}^{\ast
},y)]\leq 0$ for every $i\in I.$

Finally, we will prove that

sup$_{y\in A_{i}(\omega ,\widetilde{x}^{\ast }))_{a_{i}(\widetilde{x}^{\ast
})}}\sup_{u\in G_{i}(\omega ,\widetilde{x}^{\ast }(\omega ))_{g_{i}(%
\widetilde{x}^{\ast }(\omega ))}}[\mathit{Re}\langle u,\widetilde{x}^{\ast
}(\omega )-y\rangle +f_{i}(\omega ,\widetilde{x}^{\ast },y)]\leq 0$ for
every $i\in I.$

In order to do that, let us consider $i\in I$ and the fixed point $\omega
\in \Omega .$

Let $y\in (A_{i}(\omega ,\widetilde{x}^{\ast }))_{a_{i}(\widetilde{x}^{\ast
})}$, $\lambda \in \lbrack 0,1]$ and $z_{\lambda }^{i}(\omega ):=\lambda
y+(1-\lambda )\widetilde{x_{i}}^{\ast }(\omega ).$ According to assumption
A2 b), $z_{\lambda }^{i}(\omega )\in A_{i}(\omega ,\widetilde{x}^{\ast }).$

According to (1), we have $\sup_{u\in (G_{i}(\omega ,z_{\lambda }^{i}(\omega
)))_{g_{i}(z_{\lambda }^{i}(\omega ))}}[\mathit{Re}\langle u,\widetilde{x}%
^{\ast }(\omega )-z_{\lambda }(\omega )\rangle +f_{i}(\omega ,\widetilde{x}%
^{\ast },z_{\lambda }(\omega ))]\leq 0$ for each $\lambda \in \lbrack 0,1]$.

According to assumption A4 a), $f_{i}(\omega ,\widetilde{x}^{\ast },%
\widetilde{x}^{\ast }(t,\omega ))=0$\textit{.} For each $y_{1},y_{2}\in Y$\
and for each $\lambda \in \lbrack 0,1],$\ we also have that $f_{i}(\omega ,%
\widetilde{x}^{\ast },\lambda y_{1}+(1-\lambda )y_{2})\geq \lambda
f_{i}(\omega ,\widetilde{x}^{\ast },y_{1})+(1-\lambda )f_{i}(\omega ,%
\widetilde{x}^{\ast },y_{2}).$

Therefore, for each $\lambda \in \lbrack 0,1]$, we have that

$t\{\sup_{u\in (G_{i}(\omega ,z_{\lambda }^{i}(\omega )))_{g_{i}(z_{\lambda
}^{i}(\omega ))}}[\mathit{Re}\langle u,\widetilde{x}^{\ast }(\omega
)-y\rangle +f_{i}(\omega ,\widetilde{x}^{\ast },y)]\}\leq $

$\sup_{u\in F_{i}(\omega ,(G_{i}(\omega ,z_{\lambda }^{i}(\omega
)))_{g_{i}(z_{\lambda }^{i}(\omega ))})}t[\mathit{Re}\langle u,\widetilde{x}%
^{\ast }(\omega _{i})-y)\rangle +f_{i}(\omega ,\widetilde{x}^{\ast
},z_{\lambda }^{i}(t,\omega ))]=$

$\sup_{u\in (G_{i}(\omega ,z_{\lambda }^{i}(\omega )))_{g_{i}(z_{\lambda
}^{i}(\omega ))})}[\mathit{Re}\langle u,\widetilde{x}^{\ast }(\omega
)-z_{\lambda }^{i}(\omega )\rangle +f_{i}(\omega ,\widetilde{x}^{\ast
},z_{\lambda }^{i}(\omega ))]\leq 0.$

It follows that for each $\lambda \in \lbrack 0,1],$

(2) $\sup_{u\in (G_{i}(\omega ,z_{\lambda }^{i}(\omega )))_{g_{i}(z_{\lambda
}^{i}(\omega ))}}[\mathit{Re}\langle u,\widetilde{x}^{\ast }(\omega
)-y\rangle +f_{i}(\omega ,\widetilde{x}^{\ast },y)]\leq 0.$

Now, we are using the lower semicontinuity of $y\rightarrow (G_{i}(\omega
,y))_{g_{i}(y)}:L\cap Y\rightarrow 2^{Y^{\prime }}$ in order to show the
conclusion. For each $z_{0}\in (G_{i}(\omega ,\widetilde{x}^{\ast }(\omega
)))_{g_{i}(\widetilde{x}^{\ast }(\omega ))}$ and $e>0$ let us consider $%
U_{z_{0}}^{i},$ the neighborhood of $z_{0}$ in the topology $\sigma
(Y^{\prime },Y),$ defined by $U_{z_{0}}^{i}:=\{z\in Y^{\prime }:|\func{Re}%
\langle z_{0}-z,\widetilde{x}^{\ast }(\omega )-y\rangle |<e\}.$ As $%
y\rightarrow (G_{i}(\omega ,y))_{g_{i}(y)}:L\cap Y\rightarrow 2^{Y^{\prime
}} $ is lower semicontinuous, where $L=\{z_{\lambda }(\omega ):\lambda \in
\lbrack 0,1]\}$ and $U_{z_{0}}^{i}\cap (G_{i}(\omega ,\widetilde{x}^{\ast
}(\omega )))_{g_{i}(\widetilde{x}^{\ast }(\omega ))}\neq \emptyset ,$ there
exists a non-empty neighborhood $N(\widetilde{x}^{\ast }(\omega ))$ of $%
\widetilde{x}^{\ast }(\omega )$ in $L$ such that for each $z\in N(\widetilde{%
x}^{\ast }(\omega )),$ we have that $U_{z_{0}}^{i}\cap (G_{i}(\omega
,z))_{g_{i}(z)}\neq \emptyset .$ Then there exists $\delta \in (0,1],$ $t\in
(0,\delta )$ and $u\in (G_{i}(\omega ,z_{\lambda }^{i}(\omega
)))_{g_{i}(z_{\lambda }^{i}(\omega ))}\cap U_{z_{0}}^{i}\neq \emptyset $
such that $\mathit{Re}\langle z_{0}-u,\widetilde{x}^{\ast }(\omega
)-y\rangle <e.$ Therefore, $\mathit{Re}\langle z_{0},\widetilde{x}^{\ast
}(\omega )-y\rangle <\mathit{Re}\langle u_{i},\widetilde{x}^{\ast }(\omega
)-y\rangle +e.$

It follows that

$\mathit{Re}\langle z_{0},\widetilde{x}^{\ast }(\omega )-y\rangle +f(\omega ,%
\widetilde{x}^{\ast },y)<\mathit{Re}\langle u,\widetilde{x}^{\ast }(\omega
)-y\rangle +f(\omega ,\widetilde{x}^{\ast },y)+e<e.$

The last inequality comes from (2). Since $e>0$ and $z_{0}\in (G_{i}(\omega ,%
\widetilde{x}^{\ast }(\omega )))_{g_{i}(\widetilde{x}^{\ast }(\omega ))}$
have been chosen arbitrarily, the next relation holds:

$\mathit{Re}\langle z_{0},\widetilde{x}^{\ast }(\omega )-y\rangle
+f_{i}(\omega ,\widetilde{x}^{\ast },y)<0.$

Hence, for ecah $i\in I,$ we have that $\sup_{u\in (G_{i}(\omega ,\widetilde{%
x}^{\ast }(\omega )))_{g_{i}(\widetilde{x}^{\ast }(\omega ))}}$[$\mathit{Re}%
\langle z_{0},\widetilde{x}^{\ast }(\omega )-y\rangle +f_{i}(\omega ,%
\widetilde{x}^{\ast },y)]\leq 0$ for every $y\in $cl$(A_{i}(\omega ,%
\widetilde{x}^{\ast }))_{a_{i}(\widetilde{x}^{\ast })}$.$\medskip $

If \TEXTsymbol{\vert}I\TEXTsymbol{\vert}=1, we obtain the following
corollary.

\begin{corollary}
\textit{Let }$(\Omega $\textit{, }$\tciFourier ,\mu )$ be a complete finite
separable measure space and let $Y$ be a separable Banach space.\textit{\
Suppose that the following conditions are satisfied:}
\end{corollary}

\textit{A.1)}

\textit{\ \ \ \ \ (a) }$X:\Omega \rightarrow \mathcal{F}(Y)$\textit{\ is
such that }$\omega \rightarrow (X(\omega )_{z}:\Omega \rightarrow 2^{Y}$ 
\textit{is} \textit{a nonempty, convex, weakly compact-valued and integrably
bounded correspondence.}

\textit{\ \ \ \ \ (b) }$X:\Omega \rightarrow \mathcal{F}(Y)$\textit{\ is
such that }$\omega \rightarrow (X(\omega ))_{z}:\Omega \rightarrow 2^{Y}$ 
\textit{is }$\tciFourier -$\textit{lower measurable;}

\textit{\ A.2)}

\textit{\ \ \ \ \ (a) For each }$(\omega ,\widetilde{x})\in \Omega \times
L_{X},$\textit{\ }$(A(\omega ,\widetilde{x}))_{a(\widetilde{x})}$\textit{\
is\ convex and has a non-empty interior in the relative norm topology of }$%
(X(\omega ))_{z}.$\textit{\medskip }

\textit{\ \ \ \ \ (b)\ the correspondence }$(\omega ,\widetilde{x}%
)\rightarrow (A(\omega ,\widetilde{x}))_{a(\widetilde{x})}:\Omega \times
L_{X}\rightarrow 2^{Y}$\textit{\ has measurable graph i.e. }$\{(\omega ,%
\widetilde{x},y)\in \Omega \times L_{X}\times Y:y\in (A(\omega ,\widetilde{x}%
))_{a(\widetilde{x})}\}\in \tciFourier \otimes $\textit{\ss }$%
_{w}(L_{X})\otimes $\textit{\ss }$(Y)$\textit{\ where \ss }$_{w}(L_{X})$%
\textit{\ is the Borel }$\sigma -$\textit{algebra for the weak topology on }$%
L_{X}$\textit{\ and \ss }$(Y)$\textit{\ is the Borel }$\sigma -$\textit{%
algebra for the norm topology on }$Y$\textit{.}

\textit{\ \ \ \ \ (c)\ the correspondence }$(\omega ,\widetilde{x}%
)\rightarrow (A(\omega ,\widetilde{x}))_{a(\widetilde{x})}$\textit{\ has
weakly open lower sections, i.e., for each }$\omega \in \Omega $\textit{\
and for each }$y\in Y,$\textit{\ the set }$((A(\omega ,\widetilde{x})_{a(%
\widetilde{x})})^{-1}(\omega ,y)=\{\widetilde{x}\in L_{X}:y\in (A(\omega ,%
\widetilde{x}))_{a(\widetilde{x})}\}\}$\textit{\ is weakly open in }$L_{X};$

\textit{\ \ \ \ \ (d) For each }$\omega \in \Omega ,$\textit{\ }$\widetilde{x%
}\rightarrow $cl$(A(\omega ,\widetilde{x}))_{a(\widetilde{x}%
)}:L_{X}\rightarrow 2^{Y}$\textit{\ is upper semicontinuous in the sense
that the set }$\{\widetilde{x}\in L_{X}\mathit{:}$cl$(A(\omega ,\widetilde{x}%
))_{a(\widetilde{x})}\subset V$ \textit{is weakly open in }$L_{X}$\textit{\
for every norm open subset }$V$\textit{\ of }$Y$\textit{;}

\textit{A.3)}

$G:\Omega \times Y\rightarrow \mathcal{F}(Y^{\prime })$ \textit{and} $%
g:Y\rightarrow (0,1]$ \textit{are such that:}

\ \ \ (\textit{a) for each }$\omega \in \Omega ,$\textit{\ }$y\rightarrow
(G(\omega ,y))_{g(y)}:Y\rightarrow 2^{Y^{\prime }}$ \textit{is monotone\
with non-empty values}$;$

\ \ \ (\textit{b) for each }$\omega \in \Omega ,$\textit{\ }$y\rightarrow
(G(\omega ,y))_{g(y)}:L\cap Y\rightarrow 2^{Y^{\prime }}$\textit{\ is lower
semicontinuous from the relative topology of }$Y$\textit{\ into the weak}$%
^{\ast }-$\textit{topology }$\sigma (Y^{\prime },Y)$ \textit{of }$Y^{\prime
} $\textit{\ for each one-dimensional flat }$L\subset Y;$

\ \ \ \textit{A.4)}

\ \ \ (\textit{a) }$f:\Omega \times L_{X}\times Y\rightarrow R\cup \{\infty
,-\infty \}$\textit{\ is such that }$\widetilde{x}\rightarrow f(\omega ,%
\widetilde{x},y)$\textit{\ is lower semicontinuous on }$L_{X}$\textit{\ for
each fixed }$(\omega ,y)\in \Omega \times Y,$\textit{\ }$f(\omega ,%
\widetilde{x},\widetilde{x}(t,\omega ))=0$\textit{\ for each }$(\omega ,%
\widetilde{x})\in \Omega \times L_{X}$\textit{\ and }$y\rightarrow f(\omega ,%
\widetilde{x},y)$\textit{\ is concave on }$Y$ \textit{for each fixed }$%
(\omega ,\widetilde{x})\in \Omega \times L_{X}$\textit{;}

\ \ \ (\textit{b) for each fixed }$\omega \in \Omega ,$\textit{\ the set}

\textit{\ }$\{\widetilde{x}\in S_{X}^{1}:\sup_{y\in (A(\omega ,\widetilde{x}%
))_{a(\widetilde{x})})}[\sup_{u\in (G(\omega ,y))_{g(y)}}$\textit{Re}$%
\langle u,\widetilde{x}-y\rangle +f(\omega ,\widetilde{x},y)]>0\}$\textit{\
is weakly open in }$L_{X}$

\ \ \ (c) $\{(\omega ,\widetilde{x}):\sup_{u\in (G(\omega ,y))_{g(y)}}%
\mathit{Re}\langle u,\widetilde{x}-y\rangle +f_{i}(\omega ,\widetilde{x}%
,y)>0\}\in \mathcal{F}\otimes B(L_{X})$\textit{.}

\textit{Then, there exists }$\widetilde{x}^{\ast }\in L_{X}$\textit{\ such
that}:

i) $\widetilde{x}^{\ast }(\omega )\in $cl$(A(\omega ,\widetilde{x}^{\ast
}))_{a(\widetilde{x}^{\ast })};$

ii) \textit{sup}$_{u\in (G(\omega ,\widetilde{x}^{\ast }(\omega )))_{g(%
\widetilde{x}^{\ast }(\omega ))}}Re\langle u,\widetilde{x}^{\ast }(\omega
)-y\rangle +f(\omega ,\widetilde{x},y)]\leq 0$\textit{\ for all }$y\in
(A(\omega ,\widetilde{x}^{\ast }))_{a(\widetilde{x}^{\ast })}))$\textit{.}$%
\medskip $

We obtain the following random fixed point theorem as a particular case of
Theorem 3$.$

\begin{theorem}
\textit{Let }$(\Omega $\textit{, }$\tciFourier ,\mu )$ be a complete finite
separable measure space and let $Y$ be a separable Banach space.\textit{\
Suppose that the following conditions are satisfied:}
\end{theorem}

\textit{For each }$i\in I:$

\textit{A.1)}

\textit{\ \ \ \ \ (a) }$X_{i}:\Omega \rightarrow \mathcal{F}(Y)$\textit{\ is
such that }$\omega \rightarrow (X_{i}(\omega ))_{z}:\Omega \rightarrow 2^{Y}$
\textit{is} \textit{a nonempty, convex, weakly compact-valued and integrably
bounded correspondence.}

\textit{\ \ \ \ \ (b) }$X_{i}:\Omega \rightarrow \mathcal{F}(Y)$\textit{\ is
such that }$\omega \rightarrow (X_{i}(\omega ))_{z}:\Omega \rightarrow 2^{Y}$
\textit{is }$\tciFourier _{i}-$\textit{lower measurable;}

\textit{\ A.2)}

\textit{\ \ \ \ \ (a) For each }$(\omega ,\widetilde{x})\in \Omega \times
L_{X},$\textit{\ }$(A_{i}(\omega ,\widetilde{x}))_{a_{i}(\widetilde{x})}$%
\textit{\ is\ convex and has a non-empty interior in the relative norm
topology of }$(X_{i}(\omega ))_{z}.$\textit{\medskip }

\textit{\ \ \ \ \ (b)\ the correspondence }$(\omega ,\widetilde{x}%
)\rightarrow (A_{i}(\omega ,\widetilde{x}))_{a_{i}(\widetilde{x})}:\Omega
\times L_{X}\rightarrow 2^{Y}$\textit{\ has measurable graph i.e. }$%
\{(\omega ,\widetilde{x},y)\in \Omega \times L_{X}\times Y:y\in
(A_{i}(\omega ,\widetilde{x}))_{a_{i}(\widetilde{x})}\}\in \tciFourier
\otimes $\textit{\ss }$_{w}(L_{X})\otimes $\textit{\ss }$(Y)$\textit{\ where 
\ss }$_{w}(L_{X})$\textit{\ is the Borel }$\sigma -$\textit{algebra for the
weak topology on }$L_{X}$\textit{\ and \ss }$(Y)$\textit{\ is the Borel }$%
\sigma -$\textit{algebra for the norm topology on }$Y$\textit{.}

\textit{\ \ \ \ \ (c)\ the correspondence }$(\omega ,\widetilde{x}%
)\rightarrow (A_{i}(\omega ,\widetilde{x}))_{a_{i}(\widetilde{x})}$\textit{\
has weakly open lower sections, i.e., for each }$\omega \in \Omega $\textit{%
\ and for each }$y\in Y,$\textit{\ the set }$((A_{i}(\omega ,\widetilde{x}%
)_{a_{i}(\widetilde{x})})^{-1}(\omega ,y)=\{\widetilde{x}\in L_{X}:y\in
(A_{i}(\omega ,\widetilde{x}))_{a_{i}(\widetilde{x})}\}\}$\textit{\ is
weakly open in }$L_{X};$

\textit{\ \ \ \ \ (d) For each }$\omega \in \Omega ,$\textit{\ }$\widetilde{x%
}\rightarrow $cl$(A_{i}(\omega ,\widetilde{x}))_{a_{i}(\widetilde{x}%
)}:L_{X}\rightarrow 2^{Y}$\textit{\ is upper semicontinuous in the sense
that the set }$\{\widetilde{x}\in L_{X}\mathit{:}$cl$(A_{i}(\omega ,%
\widetilde{x}))_{a_{i}(\widetilde{x})}\subset V\}$ \textit{is weakly open in 
}$L_{X}$\textit{\ for every norm open subset }$V$\textit{\ of }$Y$\textit{;}

\textit{Then, there exists }$\widetilde{x}^{\ast }\in L_{X}$\textit{\ such
that }for every $i\in I$, $\widetilde{x_{i}}^{\ast }(\omega )\in $cl$%
(A_{i}(\omega ,\widetilde{x}^{\ast }))_{a_{i}(\widetilde{x}^{\ast
})}.\medskip $

If \TEXTsymbol{\vert}I\TEXTsymbol{\vert}=1, we obtain the following result.

\begin{theorem}
\textit{Let }$(\Omega $\textit{, }$\tciFourier ,\mu )$ be a complete finite
separable measure space and let $Y$ be a separable Banach space.\textit{\
Suppose that the following conditions are satisfied:}
\end{theorem}

\textit{\ (a) }$X:\Omega \rightarrow \mathcal{F}(Y)$\textit{\ is such that }$%
\omega \rightarrow (X(\omega ))_{z}:\Omega \rightarrow 2^{Y}$ \textit{is} 
\textit{a nonempty, convex, weakly compact-valued and integrably bounded
correspondence.}

\textit{\ \ \ \ \ (b) }$X:\Omega \rightarrow \mathcal{F}(Y)$\textit{\ is
such that }$\omega \rightarrow (X(\omega ))_{z}:\Omega \rightarrow 2^{Y}$ 
\textit{is }$\tciFourier -$\textit{lower measurable;}

\textit{\ A.2)}

\textit{\ \ \ \ \ (a) For each }$(\omega ,\widetilde{x})\in \Omega \times
L_{X},$\textit{\ }$(A(\omega ,\widetilde{x}))_{a(\widetilde{x})}$\textit{\
is\ convex and has a non-empty interior in the relative norm topology of }$%
(X(\omega ))_{z}.$\textit{\medskip }

\textit{\ \ \ \ \ (b)\ the correspondence }$(\omega ,\widetilde{x}%
)\rightarrow (A(\omega ,\widetilde{x}))_{a(\widetilde{x})}:\Omega \times
L_{X}\rightarrow 2^{Y}$\textit{\ has measurable graph i.e. }$\{(\omega ,%
\widetilde{x},y)\in \Omega \times L_{X}\times Y:y\in (A(\omega ,\widetilde{x}%
))_{a_{i}(\widetilde{x})}\}\in \tciFourier \otimes $\textit{\ss }$%
_{w}(L_{X})\otimes $\textit{\ss }$(Y)$\textit{\ where \ss }$_{w}(L_{X})$%
\textit{\ is the Borel }$\sigma -$\textit{algebra for the weak topology on }$%
L_{X}$\textit{\ and \ss }$(Y)$\textit{\ is the Borel }$\sigma -$\textit{%
algebra for the norm topology on }$Y$\textit{.}

\textit{\ \ \ \ \ (c)\ the correspondence }$(\omega ,\widetilde{x}%
)\rightarrow (A(\omega ,\widetilde{x}))_{a(\widetilde{x})}$\textit{\ has
weakly open lower sections, i.e., for each }$\omega \in \Omega $\textit{\
and for each }$y\in Y,$\textit{\ the set }$((A(\omega ,\widetilde{x})_{a(%
\widetilde{x})})^{-1}(\omega ,y)=\{\widetilde{x}\in L_{X}:y\in (A(\omega ,%
\widetilde{x}))_{a_{i}(\widetilde{x})}\}\}$\textit{\ is weakly open in }$%
L_{X};$

\textit{\ \ \ \ \ (d) For each }$\omega \in \Omega ,$\textit{\ }$\widetilde{x%
}\rightarrow $cl$(A(\omega ,\widetilde{x}))_{a(\widetilde{x}%
)}:L_{X}\rightarrow 2^{Y}$\textit{\ is upper semicontinuous in the sense
that the set }$\{\widetilde{x}\in L_{X}\mathit{:}$cl$(A(\omega ,\widetilde{x}%
))_{a(\widetilde{x})}\subset V\}$ \textit{is weakly open in }$L_{X}$\textit{%
\ for every norm open subset }$V$\textit{\ of }$Y$\textit{;}

\textit{Then, there exists }$\widetilde{x}^{\ast }\in L_{X}$\textit{\ such
that} \ \ \ \ \ \ \ \ \ \ \ \ \ \ \ \ \ \ \ \ \ \ \ \ \ \ \ \ \ \ \ \ \ $%
\widetilde{x}^{\ast }(\omega )\in $cl$(A(\omega ,\widetilde{x}^{\ast }))_{a(%
\widetilde{x}^{\ast })}.\medskip $

\begin{theorem}
\textit{Let }$I$\textit{\ be a countable or uncounatble set. Let }$(\Omega $%
\textit{, }$\tciFourier ,\mu )$ be a complete finite separable measure space
and let $Y$ be a separable Banach space.\textit{\ Suppose that the following
conditions are satisfied:}
\end{theorem}

\textit{For each }$i\in I:$

\textit{\ \ \ \ \ (a) }$X_{i}:\Omega \rightarrow \mathcal{F}(Y)$\textit{\ is
such that }$\omega \rightarrow (X_{i}(\omega ))_{z_{i}}:\Omega \rightarrow
2^{Y}$ \textit{is} \textit{a nonempty, convex, weakly compact-valued and
integrably bounded correspondence;}

\textit{\ \ \ \ \ (b) }$X_{i}:\Omega \rightarrow \mathcal{F}(Y)$\textit{\ is
such that }$\omega \rightarrow (X_{i}(\omega ))_{z_{i}}:\Omega \rightarrow
2^{Y}$ \textit{is }$\tciFourier _{i}-$\textit{lower measurable;}

\textit{\ A.2)}

\textit{\ \ \ \ \ (a) For each }$(\omega ,\widetilde{x})\in \Omega \times
L_{X},$\textit{\ }$(A_{i}(\omega ,\widetilde{x}))_{a_{i}(\widetilde{x})}$%
\textit{\ is\ nonempty convex and compact}$;$\textit{\medskip }

\textit{\ \ \ \ \ (b)\ For each }$\widetilde{x}\in L_{X},$ \textit{the
correspondence }$\omega \rightarrow (A_{i}(\omega ,\widetilde{x}))_{a_{i}(%
\widetilde{x})}:\Omega \rightarrow 2^{Y}$\textit{\ has a measurable graph;}

\textit{\ \ \ \ (c) For each }$\omega \in \Omega ,$\textit{\ }$\widetilde{x}%
\rightarrow (A_{i}(\omega ,\widetilde{x}))_{a_{i}(\widetilde{x}%
)}:L_{X}\rightarrow 2^{Y}$\textit{\ is upper semicontinuous};

\textit{A.3)}

\ \ \ \textit{\ }$\psi _{i}:\Omega \times L_{X}\times Y\rightarrow R\cup
\{-\infty ,+\infty \}$\textit{\ is such that:}

\ \ \ (\textit{a) }$\widetilde{x}\rightarrow \{y\in Y:\psi _{i}(\omega ,%
\widetilde{x},y)>0\}:L_{X}\rightarrow 2^{Y}$\textit{\ is upper
semicontinuous with compact values on }$L_{X}$\textit{\ for each fixed }$%
\omega \in \Omega ;$

\ \ \ (\textit{b) }$\widetilde{x_{i}}(\omega )\notin \{y\in Y:\psi
_{i}(\omega ,\widetilde{x},y)>0\}$\textit{\ for each fixed }$(\omega ,%
\widetilde{x})\in \Omega \times L_{X};$

\ \ \ (\textit{c) for each} $(\omega ,\widetilde{x})\in \Omega \times L_{X},$
$\psi _{i}(\omega ,\widetilde{x},\cdot )$ \textit{is concave;}

\ \ \ (\textit{d) for each }$\omega \in \Omega ,$\textit{\ }$\{\widetilde{x}%
\in L_{X}:\alpha _{i}(\omega ,\widetilde{x})>0\}$\textit{\ is weakly open in}
$L_{X},$\textit{\ where }$\alpha _{i}:\Omega \times L_{X}\rightarrow R$%
\textit{\ is defined by }$\alpha _{i}(\omega ,\widetilde{x})=\sup_{y\in
(A_{i}(\omega ,\widetilde{x}))_{a_{i}(\widetilde{x})}}\psi _{i}(\omega ,%
\widetilde{x},y)$\textit{\ for each }$(\omega ,\widetilde{x})\in \Omega
\times L_{X};$

\ \ \ (\textit{e) }$\{(\omega ,\widetilde{x}):\alpha _{i}(\omega ,\widetilde{%
x})>0\}\in \mathcal{F}_{i}\otimes B(L_{X})$\textit{.}

\textit{Then, there exists }$\widetilde{x}^{\ast }\in L_{X}$\textit{\ such
that for every} $i\in I$,

i) $\widetilde{x_{i}}^{\ast }(\omega )\in (A_{i}(\omega ,\widetilde{x}^{\ast
}))_{a_{i}(\widetilde{x}^{\ast })};$

\textit{ii) sup}$_{y\in (A_{i}(\omega ,\widetilde{x}^{\ast }))_{a_{i}(%
\widetilde{x}^{\ast })}}\psi _{i}(\omega ,\widetilde{x}^{\ast },y)\leq 0$%
\textit{.}$\medskip $

\textit{Proof.} For every $i\in I,$ let $P_{i}:\Omega \times
L_{X}\rightarrow \mathcal{F}(Y)$ and $p_{i}:L_{X}\rightarrow (0,1]$ such
that $(P_{i}(\omega ,\widetilde{x}))_{p_{i}(\widetilde{x})}=\{y\in Y:\psi
_{i}(\omega ,\widetilde{x},y)>0\}$ for each $(\omega ,\widetilde{x})\in
\Omega \times L_{X}.$

We shall show that the abstract economy $G=\{(\Omega ,\tciFourier ,\mu ),$ $%
(X_{i},\tciFourier _{i},A_{i},P_{i},a_{i},p_{i},z_{i})_{i\in I}\}$ satisfies
all the hypotheses of Theorem 1.

Suppose $\omega \in \Omega .$

According to A3 a), we have that\textit{\ }$\widetilde{x}\rightarrow
(P_{i}(\omega ,\widetilde{x}))_{p_{i}(\widetilde{x})}:L_{X}\rightarrow 2^{Y}$%
\textit{\ }is upper semicontinuous\textit{\ }with nonempty values and
according to A3 b), $\widetilde{x_{i}}(\omega )\not\in (P_{i}(\omega ,%
\widetilde{x}))_{p_{i}(\widetilde{x})}$ for each $\widetilde{x}\in L_{X}.$
Assumption A3 c) implies that $\widetilde{x}\rightarrow (P_{i}(\omega ,%
\widetilde{x}))_{p_{i}(\widetilde{x})}:L_{X}\rightarrow 2^{Y}$ has convex
values.

By the definition of $\alpha _{i},$ we note that $\{\widetilde{x}\in
L_{X}:(A_{i}(\omega ,\widetilde{x}))_{a_{i}(\widetilde{x})}\cap
(P_{i}(\omega ,\widetilde{x}))_{p_{i}(\widetilde{x})}\neq \emptyset \}=\{%
\widetilde{x}\in L_{X}:\alpha _{i}(\omega ,\widetilde{x})>0\}$ so that $\{%
\widetilde{x}\in L_{X}:(A_{i}(\omega ,\widetilde{x}))_{a_{i}(\widetilde{x}%
)}\cap (P_{i}(\omega ,\widetilde{x}))_{p_{i}(\widetilde{x})}\neq \emptyset
\} $ is weakly open in $L_{X}$ by A3 d).

According to A2 b) and A3 e), it follows that the correspondences\textit{\ }$%
(\omega ,\widetilde{x})\rightarrow (A_{i}(\omega ,\widetilde{x}))_{a_{i}(%
\widetilde{x})}:\Omega \times L_{X}\rightarrow 2^{Y}$ and $(\omega ,%
\widetilde{x})\rightarrow (P_{i}(\omega ,\widetilde{x}))_{p_{i}(\widetilde{x}%
)}:\Omega \times L_{X}\rightarrow 2^{Y}$\textit{\ }have measurable graphs$.$

Thus the Bayesian abstract fuzzy economy $G=\{(\Omega ,\tciFourier ,\mu ),$ $%
(X_{i},\tciFourier _{i},A_{i},P_{i},a_{i},b_{i},z_{i})_{i\in I}\}$ satisfies
all the hypotheses of Theorem 2. Therefore, there exists $\widetilde{x}%
^{\ast }\in L_{X}$ such that for every $i\in I:$

$\widetilde{x}_{i}^{\ast }(\omega )\in (A_{i}(\omega ,\widetilde{x}^{\ast
}))_{a_{i}(\widetilde{x}^{\ast })}$ $\mu -a.e$ and

$(A_{i}(\omega ,\widetilde{x}^{\ast }))_{a_{i}(\widetilde{x}^{\ast })}\cap
(P_{i}(\omega ,\widetilde{x}^{\ast }))_{p_{i}(\widetilde{x}^{\ast })}=\phi $ 
$\mu -a.e;$

that is, there exists $\widetilde{x}^{\ast }\in L_{X}$ such that for every $%
i\in I:$

i) $\widetilde{x_{i}}^{\ast }(\omega )\in (A_{i}(\omega ,\widetilde{x}^{\ast
}))_{a_{i}(\widetilde{x}^{\ast })};$

ii) sup$_{y\in (A_{i}(\omega ,\widetilde{x}^{\ast }))_{a_{i}(\widetilde{x}%
^{\ast })}}\psi _{i}(\omega ,\widetilde{x}^{\ast },y)\leq 0$.$\medskip $

If \TEXTsymbol{\vert}I\TEXTsymbol{\vert}=1, we obtain the following
corollary.

\begin{corollary}
\textit{Let }$(\Omega $\textit{, }$\tciFourier ,\mu )$ be a complete finite
separable measure space and let $Y$ be a separable Banach space.\textit{\
Suppose that the following conditions are satisfied:}
\end{corollary}

\textit{A.1)}

\textit{\ \ \ \ \ (a) }$X:\Omega \rightarrow \mathcal{F}(Y)$\textit{\ is
such that }$\omega \rightarrow (X(\omega ))_{z}:\Omega \rightarrow 2^{Y}$ 
\textit{is} \textit{a nonempty, convex, weakly compact-valued and integrably
bounded correspondence.}

\textit{\ \ \ \ \ (b) }$X:\Omega \rightarrow \mathcal{F}(Y)$\textit{\ is
such that }$\omega \rightarrow (X(\omega ))_{z}:\Omega \rightarrow 2^{Y}$ 
\textit{is }$\tciFourier -$\textit{lower measurable;}

\textit{\ A.2)}

\textit{\ \ \ \ \ (a) For each }$(\omega ,\widetilde{x})\in \Omega \times
L_{X},$\textit{\ }$(A(\omega ,\widetilde{x}))_{a(\widetilde{x})}$\textit{\
is\ non-empty convex and compact}$;$\textit{\medskip }

\textit{\ \ \ \ \ (b)\ For each }$\widetilde{x}\in L_{X},$ \textit{the
correspondence }$\omega \rightarrow (A(\omega ,\widetilde{x}))_{a(\widetilde{%
x})}:\Omega \rightarrow 2^{Y}$\textit{\ has a measurable graph;}

\textit{\ \ \ \ (c) For each }$\omega \in \Omega ,$\textit{\ }$\widetilde{x}%
\rightarrow (A(\omega ,\widetilde{x}))_{a(\widetilde{x})}:L_{X}\rightarrow
2^{Y}$\textit{\ is upper semicontinuous;}

\textit{A.3)}

\ \ \ \textit{\ }$\psi :\Omega \times L_{X}\times Y\rightarrow R\cup
\{-\infty ,+\infty \}$\textit{\ is such that:}

\ \ $\ $(a) $\widetilde{x}\rightarrow \{y\in Y:\psi (\omega ,\widetilde{x}%
,y)>0\}:L_{X}\rightarrow 2^{Y}$\textit{\ is upper semicontinuous with weakly
compact values on }$L_{X}$\textit{\ for each fixed }$\omega \in \Omega ;$

\ \ \ (\textit{b) }$\widetilde{x}(\omega )\notin \{y\in Y:\psi (\omega ,%
\widetilde{x},y)>0\}$\textit{\ for each fixed }$(\omega ,\widetilde{x})\in
\Omega \times L_{X};$

\ \ \ (\textit{c) for each} $(\omega ,\widetilde{x})\in \Omega \times L_{X},$
$\psi (\omega ,\widetilde{x},\cdot )$ \textit{is concave;}

\ \ \ (\textit{d) for each }$\omega \in \Omega ,$\textit{\ }$\{\widetilde{x}%
\in L_{X}:\alpha (\omega ,\widetilde{x})>0\}$\textit{\ is weakly open in} $%
L_{X},$\textit{\ where }$\alpha :\Omega \times L_{X}\rightarrow R$\textit{\
is defined by }$\alpha (\omega ,\widetilde{x})=\sup_{y\in (A(\omega ,%
\widetilde{x}))_{a(\widetilde{x})}}\psi (\omega ,\widetilde{x},y)$\textit{\
for each }$(\omega ,\widetilde{x})\in \Omega \times L_{X};$

\ \ \ (\textit{e) }$\{(\omega ,\widetilde{x}):\alpha (\omega ,\widetilde{x}%
)>0\}\in \mathcal{F}\otimes B(L_{X})$\textit{.}

\textit{Then, there exists }$\widetilde{x}^{\ast }\in L_{X}$\textit{\ such
that}:

i) $\widetilde{x}^{\ast }(\omega )\in (A(\omega ,\widetilde{x}^{\ast }))_{a(%
\widetilde{x}^{\ast })};$

\textit{ii) sup}$_{y\in (A(\omega ,\widetilde{x}^{\ast }))_{a(\widetilde{x}%
^{\ast })}}\psi (\omega ,\widetilde{x}^{\ast },y)\leq 0$\textit{.}$\medskip $

As a consequence of the Theorem 2, we prove the following Tan and Yuan's
type (1995) random quasi-variational inequality with random fuzzy mappings.

\begin{theorem}
\textit{Let }$(\Omega $\textit{, }$\tciFourier ,\mu )$ be a complete finite
separable measure space and let $Y$ be a separable Banach space.\textit{\
Suppose that the following conditions are satisfied:}
\end{theorem}

\textit{For each }$i\in I:$

\textit{A.1)}

\textit{\ \ \ \ \ (a) }$X_{i}:\Omega \rightarrow \mathcal{F}(Y)$\textit{\ is
such that }$\omega \rightarrow (X_{i}(\omega ))_{z_{i}}:\Omega \rightarrow
2^{Y}$ \textit{is} \textit{a nonempty, convex, weakly compact-valued and
integrably bounded correspondence;}

\textit{\ \ \ \ \ (b) }$X_{i}:\Omega \rightarrow \mathcal{F}(Y)$\textit{\ is
such that }$\omega \rightarrow (X_{i}(\omega ))_{z_{i}}:\Omega \rightarrow
2^{Y}$ \textit{is }$\tciFourier _{i}-$\textit{lower measurable;}

\textit{\ A.2)}

\textit{\ \ \ \ \ (a) For each }$(\omega ,\widetilde{x})\in \Omega \times
L_{X},$\textit{\ }$(A_{i}(\omega ,\widetilde{x}))_{a_{i}(\widetilde{x})}$%
\textit{\ is\ nonempty convex and weakly compact}$;$\textit{\medskip }

\textit{\ \ \ \ \ (b)\ For each }$\widetilde{x}\in L_{X},$ \textit{the
correspondence }$\omega \rightarrow (A_{i}(\omega ,\widetilde{x}))_{a_{i}(%
\widetilde{x})}:\Omega \rightarrow 2^{Y}$\textit{\ has a measurable graph;}

\textit{\ \ \ \ (c) For each }$\omega \in \Omega ,$\textit{\ }$\widetilde{x}%
\rightarrow (A_{i}(\omega ,\widetilde{x}))_{a_{i}(\widetilde{x}%
)}:L_{X}\rightarrow 2^{Y}$\textit{\ is upper semicontinuous};

\textit{A.3)}

$G_{i}:\Omega \times Y\rightarrow \mathcal{F}(Y^{\prime })$ \textit{and} $%
g_{i}:Y\rightarrow (0,1]$ \textit{are such that:}

\ \ (a) \textit{For each fixed }$(\omega ,y)\in \Omega \times Y,$\textit{\ }$%
\widetilde{x}\rightarrow \{y\in Y:\sup_{u\in (G_{i}(\omega ,y))_{g_{i}(y)}}$%
Re$\langle u,\widetilde{x_{i}}(\omega )-y\rangle >0\}:L_{X}\rightarrow 2^{Y}$%
\textit{\ is upper semicontinuous with compact values;}

\ \ \ (\textit{b) for each fixed }$\omega \in \Omega ,$\textit{\ the set}

\textit{\ }$\{\widetilde{x}\in L_{X}:\sup_{y\in (A_{i}(\omega ,\widetilde{x}%
))_{a_{i}(\widetilde{x})})}\sup_{u\in (G_{i}(\omega ,y))_{g_{i}(y)}}$\textit{%
Re}$\langle u,\widetilde{x_{i}}(\omega )-y\rangle >0\}$\textit{\ is weakly
open in }$L_{X}$

\ \ \ (c) $\{(\omega ,\widetilde{x}):\sup_{u\in (G_{i}(\omega
,y))_{g_{i}(y)}}\mathit{Re}\langle u,\widetilde{x_{i}}(\omega )-y\rangle
>0\}\in \mathcal{F}\otimes B(L_{X})$\textit{.}

\textit{A.4)}

$H_{i}:\Omega \times Y\rightarrow \mathcal{F}(Y^{\prime })$ \textit{and} $%
h_{i}:Y\rightarrow (0,1]$ \textit{are such that:}

\ \ \textit{(a)} \textit{For each fixed }$(\omega ,y)\in \Omega \times Y,$ $%
(H_{i}(\omega ,y))_{h_{i}(y)}\subset (G_{i}(\omega ,y))_{g_{i}(y)};$

\ \ \ (\textit{b) for each }$\omega \in \Omega ,$\textit{\ }$y\rightarrow
(H_{i}(\omega ,y))_{h_{i}(y)}:Y\rightarrow 2^{Y^{\prime }}$ \textit{is
monotone\ with non-empty values}$;$

\ \ \ (\textit{c) for each }$\omega \in \Omega ,$\textit{\ }$y\rightarrow
(H_{i}(\omega ,y))_{h_{i}(y)}:L\cap Y\rightarrow 2^{Y^{\prime }}$\textit{\
is lower semicontinuous from the relative topology of }$Y$\textit{\ into the
weak}$^{\ast }-$\textit{topology }$\sigma (Y^{\prime },Y)$ \textit{of }$%
Y^{\prime }$\textit{\ for each one-dimensional flat }$L\subset Y.$

\textit{Then, there exists }$\widetilde{x}^{\ast }\in L_{X}$\textit{\ such
that for every} $i\in I$:

i) $\widetilde{x_{i}}^{\ast }(\omega )\in (A_{i}(\omega ,\widetilde{x}^{\ast
}))_{a_{i}(\widetilde{x}^{\ast })};$

ii) \textit{sup}$_{u\in (H_{i}(\omega ,\widetilde{x}^{\ast }(\omega
)))_{h_{i}(\widetilde{x}^{\ast }(\omega ))}}Re\langle u,\widetilde{x_{i}}%
^{\ast }(\omega )-y\rangle \leq 0$\textit{\ for all }$y\in (A_{i}(\omega ,%
\widetilde{x}^{\ast }))_{a_{i}(\widetilde{x}^{\ast })}$\textit{.}$\medskip $

\textit{Proof.} Let us define $\psi _{i}:\Omega \times L_{X}\times
Y\rightarrow R\cup \{-\infty ,+\infty \}$ by

$\psi _{i}(\omega ,\widetilde{x},y)=\sup_{u\in (G_{i}(\omega
,y))_{g_{i}(y)}} $\textit{Re}$\langle u,\widetilde{x_{i}}(\omega )-y\rangle $
for each $(\omega ,\widetilde{x},y)\in \Omega \times L_{X}\times Y.$

We have that $\widetilde{x_{i}}(\omega )\notin \{y\in Y:\psi _{i}(\omega ,%
\widetilde{x},y)>0\}$\textit{\ }for each fixed\textit{\ }$(\omega ,%
\widetilde{x})\in \Omega \times L_{X}$ and, as a consequence of assumption
A3 b), it follows that for each $(\omega ,\widetilde{x})\in \Omega \times
L_{X},$ $\psi _{i}(\omega ,\widetilde{x},\cdot )$ is concave$.$

All the hypotheses of Theorem 2 are satisfied. According to Theorem 2, there
exists $\widetilde{x}^{\ast }\in L_{X}$ such that $\widetilde{x_{i}}^{\ast
}(\omega )\in $cl$(A_{i}(\omega ,\widetilde{x}^{\ast }))_{a_{i}(\widetilde{x}%
^{\ast })}$ for every $i\in I.$

and

(1) \ \ sup$_{y\in A_{i}(\omega ,\widetilde{x}^{\ast }))_{a_{i}(\widetilde{x}%
^{\ast })}}\sup_{u\in (G_{i}(\omega ,y))_{g_{i}(y)}}\mathit{Re}\langle u,%
\widetilde{x_{i}}^{\ast }(\omega )-y\rangle \leq 0$ for every $i\in I.$

Finally, we will prove that

sup$_{y\in A_{i}(\omega ,\widetilde{x}^{\ast }))_{a_{i}(\widetilde{x}^{\ast
})}}\sup_{u\in H_{i}(\omega ,\widetilde{x}^{\ast }(\omega ))_{h_{i}(%
\widetilde{x}^{\ast }(\omega ))}}\mathit{Re}\langle u,\widetilde{x_{i}}%
^{\ast }(\omega )-y\rangle \leq 0$ for every $i\in I.$

In order to do that, let us consider $i\in I$ and the fixed point $\omega
\in \Omega .$

Let $y\in (A_{i}(\omega ,\widetilde{x}^{\ast }))_{a_{i}(\widetilde{x}^{\ast
})}$, $\lambda \in \lbrack 0,1]$ and $z_{\lambda }^{i}(\omega ):=\lambda
y+(1-\lambda )\widetilde{x_{i}}^{\ast }(\omega ).$ According to assumption
A2 a), $z_{\lambda }^{i}(\omega )\in A_{i}(\omega ,\widetilde{x}^{\ast }).$

According to (1), we have $\sup_{u\in (H_{i}(\omega ,z_{\lambda }^{i}(\omega
)))_{h_{i}(z_{\lambda }^{i}(\omega ))}}\mathit{Re}\langle u,\widetilde{x_{i}}%
^{\ast }(\omega )-z_{\lambda }^{i}(\omega )\rangle \leq 0$ for each $\lambda
\in \lbrack 0,1]$.

For each $\lambda \in \lbrack 0,1]$, we have that

$t\{\sup_{u\in (H_{i}(\omega ,z_{\lambda }^{i}(\omega )))_{h_{i}(z_{\lambda
}^{i}(\omega ))}}\mathit{Re}\langle u,\widetilde{x_{i}}^{\ast }(\omega
)-y\rangle \}=$

$\sup_{u\in H_{i}(\omega ,z_{\lambda }^{i}(\omega )))_{h_{i}(z_{\lambda
}^{i}(\omega ))})}t\mathit{Re}\langle u,\widetilde{x_{i}}^{\ast }(\omega
)-y)\rangle =$

$\sup_{u\in (H_{i}(\omega ,z_{\lambda }^{i}(\omega )))_{h_{i}(z_{\lambda
}^{i}(\omega ))})}\mathit{Re}\langle u,\widetilde{x_{i}}^{\ast }(\omega
)-z_{\lambda }^{i}(\omega )\rangle \leq 0.$

It follows that for each $\lambda \in \lbrack 0,1],$

(2) $\sup_{u\in (H_{i}(\omega ,z_{\lambda }^{i}(\omega )))_{h_{i}(z_{\lambda
}^{i}(\omega ))}}\mathit{Re}\langle u,\widetilde{x_{i}}^{\ast }(\omega
)-y\rangle \leq 0.$

Now, we are using the lower semicontinuity of $y\rightarrow (H_{i}(\omega
,y))_{h_{i}(y)}:L\cap Y\rightarrow 2^{Y^{\prime }}$ in order to show the
conclusion. For each $z_{0}\in (H_{i}(\omega ,\widetilde{x_{i}}^{\ast
}(\omega )))_{h_{i}(\widetilde{x}^{\ast }(\omega ))}$ and $e>0$ let us
consider $U_{z_{0}}^{i},$ the neighborhood of $z_{0}$ in the topology $%
\sigma (Y^{\prime },Y),$ defined by $U_{z_{0}}^{i}:=\{z\in Y^{\prime }:|%
\func{Re}\langle z_{0}-z,\widetilde{x_{i}}^{\ast }(\omega )-y\rangle |<e\}.$
As $y\rightarrow (H_{i}(\omega ,y))_{h_{i}(y)}:L\cap Y\rightarrow
2^{Y^{\prime }}$ is lower semicontinuous, where $L=\{z_{\lambda }^{i}(\omega
):\lambda \in \lbrack 0,1]\}$ and $U_{z_{0}}^{i}\cap (H_{i}(\omega ,%
\widetilde{x}_{i}^{\ast }(\omega )))_{h_{i}(\widetilde{x}_{i}^{\ast }(\omega
))}\neq \emptyset ,$ there exists a non-empty neighborhood $N(\widetilde{%
x_{i}}^{\ast }(\omega ))$ of $\widetilde{x_{i}}^{\ast }(\omega )$ in $L$
such that for each $z\in N(\widetilde{x_{i}}^{\ast }(\omega )),$ we have
that $U_{z_{0}}^{i}\cap (H_{i}(\omega ,z))_{h_{i}(z)}\neq \emptyset .$ Then
there exists $\delta \in (0,1],$ $t\in (0,\delta )$ and $u\in (H_{i}(\omega
,z_{\lambda }^{i}(\omega )))_{h_{i}(z_{\lambda }^{i}(\omega ))}\cap
U_{z_{0}}^{i}\neq \emptyset $ such that $\mathit{Re}\langle z_{0}-u,%
\widetilde{x_{i}}^{\ast }(\omega )-y\rangle <e.$ Therefore, $\mathit{Re}%
\langle z_{0},\widetilde{x_{i}}^{\ast }(\omega )-y\rangle <\mathit{Re}%
\langle u_{i},\widetilde{x_{i}}^{\ast }(\omega )-y\rangle +e.$

It follows that

$\mathit{Re}\langle z_{0},\widetilde{x_{i}}^{\ast }(\omega )-y\rangle <%
\mathit{Re}\langle u,\widetilde{x_{i}}^{\ast }(\omega )-y\rangle +e<e.$

The last inequality comes from (2). Since $e>0$ and $z_{0}\in (H_{i}(\omega ,%
\widetilde{x}^{\ast }(\omega )))_{h_{i}(\widetilde{x}^{\ast }(\omega ))}$
have been chosen arbitrarily, the next relation holds:

$\mathit{Re}\langle z_{0},\widetilde{x_{i}}^{\ast }(\omega )-y\rangle <0.$

Hence, for each $i\in I,$ we have that $\sup_{u\in (H_{i}(\omega ,\widetilde{%
x}^{\ast }(\omega )))_{h_{i}(\widetilde{x}^{\ast }(\omega ))}}\mathit{Re}%
\langle z_{0},\widetilde{x_{i}}^{\ast }(\omega )-y\rangle \leq 0$ for every $%
y\in (A_{i}(\omega ,\widetilde{x}^{\ast }))_{a_{i}(\widetilde{x}^{\ast })}$.$%
\medskip $

If \TEXTsymbol{\vert}I\TEXTsymbol{\vert}=1, we obtain the following
corollary.

\begin{corollary}
\textit{Let }$(\Omega $\textit{, }$\tciFourier ,\mu )$ be a complete finite
separable measure space and let $Y$ be a separable Banach space.\textit{\
Suppose that the following conditions are satisfied:}
\end{corollary}

\textit{A.1)}

\textit{\ \ \ \ \ (a) }$X:\Omega \rightarrow \mathcal{F}(Y)$\textit{\ is
such that }$\omega \rightarrow (X(\omega )_{z}:\Omega \rightarrow 2^{Y}$ 
\textit{is} \textit{a nonempty, convex, weakly compact-valued and integrably
bounded correspondence.}

\textit{\ \ \ \ \ (b) }$X:\Omega \rightarrow \mathcal{F}(Y)$\textit{\ is
such that }$\omega \rightarrow (X(\omega ))_{z}:\Omega \rightarrow 2^{Y}$ 
\textit{is }$\tciFourier -$\textit{lower measurable;}

\textit{\ A.2)}

\textit{\ \ \ \ \ (a) For each }$(\omega ,\widetilde{x})\in \Omega \times
L_{X},$\textit{\ }$(A(\omega ,\widetilde{x}))_{a(\widetilde{x})}$\textit{\
is\ non-empty convex and compact}$.$\textit{\medskip }

\textit{\ \ \ \ \ (b)\ For each }$\widetilde{x}\in L_{X},$ \textit{the
correspondence }$\omega \rightarrow (A(\omega ,\widetilde{x}))_{a(\widetilde{%
x})}:\Omega \rightarrow 2^{Y}$\textit{\ has a measurable graph;}

\textit{\ \ \ \ (c) For each }$\omega \in \Omega ,$\textit{\ }$\widetilde{x}%
\rightarrow (A(\omega ,\widetilde{x}))_{a(\widetilde{x})}:L_{X}\rightarrow
2^{Y}$\textit{\ is upper semicontinuous;}

\textit{A.3)}

$G:\Omega \times Y\rightarrow \mathcal{F}(Y^{\prime })$ \textit{and} $%
g:Y\rightarrow (0,1]$ \textit{are such that:}

\ \ (a) \textit{For each fixed }$(\omega ,y)\in \Omega \times Y,$\textit{\ }$%
\widetilde{x}\rightarrow \{y\in Y:\sup_{u\in (G(\omega ,y))_{g(y)}}Re\langle
u,\widetilde{x}(\omega )-y\rangle >0\}:L_{X}\rightarrow 2^{Y}$\textit{\ is
upper semicontinuous with compact values;}

\ \ \ (\textit{b) for each fixed }$\omega \in \Omega ,$\textit{\ the set}

\textit{\ }$\{\widetilde{x}\in L_{X}:\sup_{y\in (A(\omega ,\widetilde{x}%
))_{a(\widetilde{x})})}\sup_{u\in (G(\omega ,y))_{g(y)}}$\textit{Re}$\langle
u,\widetilde{x}(\omega )-y\rangle >0\}$\textit{\ is weakly open in }$L_{X};$

\ \ \ (c) $\{(\omega ,\widetilde{x}):\sup_{u\in (G(\omega ,y))_{g(y)}}%
\mathit{Re}\langle u,\widetilde{x}(\omega )-y\rangle >0\}\in \mathcal{F}%
\otimes B(L_{X});$

\textit{A.4)}

$H:\Omega \times Y\rightarrow \mathcal{F}(Y^{\prime })$ \textit{and} $%
h:Y\rightarrow (0,1]$ \textit{are such that:}

\ \ \textit{(a)} \textit{For each fixed }$(\omega ,y)\in \Omega \times Y,$ $%
(H(\omega ,y))_{h(y)}\subset (G(\omega ,y))_{g(y)};$

\ \ \ (\textit{b) for each }$\omega \in \Omega ,$\textit{\ }$y\rightarrow
(H(\omega ,y))_{h(y)}:Y\rightarrow 2^{Y^{\prime }}$ \textit{is monotone\
with non-empty values}$;$

\ \ \ (\textit{c) for each }$\omega \in \Omega ,$\textit{\ }$y\rightarrow
(H(\omega ,y))_{h(y)}:L\cap Y\rightarrow 2^{Y^{\prime }}$\textit{\ is lower
semicontinuous from the relative topology of }$Y$\textit{\ into the weak}$%
^{\ast }-$\textit{topology }$\sigma (Y^{\prime },Y)$ \textit{of }$Y^{\prime
} $\textit{\ for each one-dimensional flat }$L\subset Y.$

\textit{Then, there exists }$\widetilde{x}^{\ast }\in L_{X}$\textit{\ such
that}:

i) $\widetilde{x}^{\ast }(\omega )\in (A(\omega ,\widetilde{x}^{\ast }))_{a(%
\widetilde{x}^{\ast })};$

ii) \textit{sup}$_{u\in (H(\omega ,\widetilde{x}^{\ast }(\omega )))_{h(%
\widetilde{x}^{\ast }(\omega ))}}Re\langle u,\widetilde{x}^{\ast }(\omega
)-y\rangle \leq 0$\textit{\ for all }$y\in (A(\omega ,\widetilde{x}^{\ast
}))_{a(\widetilde{x}^{\ast })}$\textit{.}$\medskip $

We obtain the following random fixed point theorem as a corollary$.$

\begin{theorem}
\textit{Let }$(\Omega $\textit{, }$\tciFourier ,\mu )$ be a complete finite
separable measure space and let $Y$ be a separable Banach space.\textit{\
Suppose that the following conditions are satisfied:}
\end{theorem}

\textit{For each }$i\in I:$

\textit{A.1)}

\textit{\ \ \ \ \ (a) }$X_{i}:\Omega \rightarrow \mathcal{F}(Y)$\textit{\ is
such that }$\omega \rightarrow (X_{i}(\omega ))_{z_{i}}:\Omega \rightarrow
2^{Y}$ \textit{is} \textit{a nonempty, convex, weakly compact-valued and
integrably bounded correspondence;}

\textit{\ \ \ \ \ (b) }$X_{i}:\Omega \rightarrow \mathcal{F}(Y)$\textit{\ is
such that }$\omega \rightarrow (X_{i}(\omega ))_{z_{i}}:\Omega \rightarrow
2^{Y}$ \textit{is }$\tciFourier _{i}-$\textit{lower measurable;}

\textit{\ A.2)}

\textit{\ \ \ \ \ (a) For each }$(\omega ,\widetilde{x})\in \Omega \times
L_{X},$\textit{\ }$(A_{i}(\omega ,\widetilde{x}))_{a_{i}(\widetilde{x})}$%
\textit{\ is\ non-empty convex and compact;\medskip }

\textit{\ \ \ \ \ (b)\ For each }$\widetilde{x}\in L_{X},$ \textit{the
correspondence }$\omega \rightarrow (A_{i}(\omega ,\widetilde{x}))_{a_{i}(%
\widetilde{x})}:\Omega \rightarrow 2^{Y}$\textit{\ has a measurable graph;}

\textit{\ \ \ \ (c) For each }$\omega \in \Omega ,$\textit{\ }$\widetilde{x}%
\rightarrow (A_{i}(\omega ,\widetilde{x}))_{a_{i}(\widetilde{x}%
)}:L_{X}\rightarrow 2^{Y}$\textit{\ is upper semicontinuous};

\textit{Then, there exists }$\widetilde{x}^{\ast }\in L_{X}$\textit{\ such
that for every} $i\in I$, $\widetilde{x_{i}}^{\ast }(\omega )\in
(A_{i}(\omega ,\widetilde{x}^{\ast }))_{a_{i}(\widetilde{x}^{\ast
})}.\medskip $

If \TEXTsymbol{\vert}I\TEXTsymbol{\vert}=1, we obtain the following result.

\begin{theorem}
\textit{Let }$(\Omega $\textit{, }$\tciFourier ,\mu )$ be a complete finite
separable measure space and let $Y$ be a separable Banach space.\textit{\
Suppose that the following conditions are satisfied:}
\end{theorem}

\textit{\ A.1)}

\textit{\ \ \ \ \ (a) }$X:\Omega \rightarrow \mathcal{F}(Y)$\textit{\ is
such that }$\omega \rightarrow (X(\omega ))_{z}:\Omega \rightarrow 2^{Y}$ 
\textit{is} \textit{a nonempty, convex, weakly compact-valued and integrably
bounded correspondence;}

\textit{\ \ \ \ \ (b) }$X:\Omega \rightarrow \mathcal{F}(Y)$\textit{\ is
such that }$\omega \rightarrow (X(\omega ))_{z}:\Omega \rightarrow 2^{Y}$ 
\textit{is }$\tciFourier -$\textit{lower measurable;}

\textit{\ A.2)}

\textit{\ \ \ \ \ (a) For each }$(\omega ,\widetilde{x})\in \Omega \times
L_{X},$\textit{\ }$(A(\omega ,\widetilde{x}))_{a(\widetilde{x})}$\textit{\
is\ non-empty convex and compact}$;$\textit{\medskip }

\textit{\ \ \ \ \ (b)\ For each }$\widetilde{x}\in L_{X},$ \textit{the
correspondence }$\omega \rightarrow (A(\omega ,\widetilde{x}))_{a(\widetilde{%
x})}:\Omega \rightarrow 2^{Y}$\textit{\ has a measurable graph;}

\textit{\ \ \ \ (c) For each }$\omega \in \Omega ,$\textit{\ }$\widetilde{x}%
\rightarrow (A(\omega ,\widetilde{x}))_{a(\widetilde{x})}:L_{X}\rightarrow
2^{Y}$\textit{\ is upper semicontinuous};

\textit{Then, there exists }$\widetilde{x}^{\ast }\in L_{X}$\textit{\ such
that for every} $i\in I$, $\widetilde{x}^{\ast }(\omega )\in (A(\omega ,%
\widetilde{x}^{\ast }))_{a(\widetilde{x}^{\ast })}\medskip $

\section{APPENDIX}

The results below have been used in the proof of our theorems. For more
details and further references see the paper quoted.

\begin{theorem}
\textit{(Projection theorem)}\textbf{.} \textit{Let }$(\Omega ,$\textit{\ }$%
\tciFourier $\textit{, }$\mu )$\textit{\ be a complete, finite measure
space, and }$Y$\textit{\ be a complete separable metric space. If }$H$%
\textit{\ belongs to }$\tciFourier \otimes $\textit{\ss }$(Y)$\textit{, its
projection Proj}$_{\Omega }(H)$\textit{\ belongs to }$\tciFourier .$
\end{theorem}

\begin{theorem}
\textit{(Aumann measurable selection theorem} [30]). \textit{Let }$(\Omega $%
\textit{, }$\tciFourier ,\mu )$\textit{\ be a complete finite measure space, 
}$Y$\textit{\ be a complete, separable metric space and }$T:\Omega
\rightarrow 2^{Y}$\textit{\ be a nonempty valued correspondence with a
measurable graph, i.e., }$G_{T}\in \tciFourier \otimes \beta (Y).$\textit{\
Then there is a measurable function }$f:\Omega \rightarrow Y$\textit{\ such
that }$f(\omega )\in T(\omega )$\textit{\ }$\mu -a.e.\medskip $
\end{theorem}

\begin{theorem}
\textit{(Diestel's Theorem [30, Theorem 3.1)}.\textit{\ Let }$(\Omega $%
\textit{, }$\tciFourier ,\mu )$\textit{\ be a complete finite measure space, 
}$X$\textit{\ be a separable Banach space and }$T:\Omega \rightarrow 2^{Y}$%
\textit{\ be an integrably bounded, convex, weakly compact and nonempty
valued correspondence. Then }$S_{T}=\{x\in L_{1}(\mu ,Y):$ $x(\omega )\in
T(\omega )$ $\mu $-a.e.\}\textit{\ is weakly compact in }$L_{1}(\mu ,Y).$%
\textit{\smallskip }
\end{theorem}

\begin{theorem}
(Carath\'{e}odory-type selection theorem [15]\textbf{)}. \textit{Let }$%
(\Omega $\textit{,}$\tciFourier ,\mu )$\textit{\ be a complete measure
space, }$Z$\textit{\ be a complete separable metric space and }$Y$\textit{\
a separable Banach space. Let }$X:\Omega \rightarrow 2^{Y}$\textit{\ be a
correspondence with a measurable graph, i.e., }$G_{X}\in \tciFourier \otimes 
$\textit{\ss }$(Y)$\textit{\ and let }$T:\Omega \times Z\rightarrow 2^{Y}$%
\textit{\ be a convex valued correspondence (possibly empty) with a
meaurable graph, i.e., }$G_{T}\in \tciFourier \otimes $\textit{\ss }$%
(Z)\otimes $\textit{\ss }$(Y)$\textit{\ where }\text{\text{\ss }}$(Y)$%
\textit{\ and \ss }$(Z)$\textit{\ are the Borel }$\sigma -$\textit{algebras
of }$Y$\textit{\ and }$Z$\textit{, respectively}$.$
\end{theorem}

\textit{Suppose that:}

\textit{(a) for each }$\omega \in \Omega $\textit{, }$T(\omega ,x)\subset
X(\omega )$\textit{\ for all }$x\in Z.$

\textit{(b) for each }$\omega \in \Omega $\textit{, }$T(\omega ,$\textit{%
\textperiodcentered }$)$\textit{\ has open lower sections in Z, i.e., for
each }$\omega \in \Omega $\textit{\ and }$y\in Y$\textit{, }$T^{-1}(\omega
,y)=\{x\in Z:y\in T(\omega ,x)\}$\textit{\ is open in Z.}

\textit{(c) for each }$(\omega ,x)\in \Omega \times Z,$\textit{\ if }$%
T(\omega ,x)\neq \emptyset $\textit{, then }$T(\omega ,x)$\textit{\ has a
non-empty interior in }$X(\omega )$\textit{.}

\textit{Let }$U=\{(\omega ,x)\in \Omega \times Z:T(\omega ,x)\neq \emptyset
\}$\textit{\ and for each }$x\in Z$\textit{, }$U^{x}=\{\omega \in \Omega
:(\omega ,x)\in U\}$\textit{\ and for each }$\omega \in \Omega $\textit{, }$%
U^{\omega }=\{x\in Z:(\omega ,x)\in U\}.$\textit{\ Then for each }$x\in Z,$%
\textit{\ }$U^{x}$\textit{\ is a measurable set in }$\Omega $\textit{\ and
there exists a Caratheodory-type selection from }$T_{\mid U}$, \textit{\
i.e., there exists a function }$f:U\rightarrow Y$\textit{\ such that }$%
f(\omega ,x)\in T(\omega ,x)$\textit{\ for all }$(\omega ,x)\in U$\textit{\
, for each }$x\in Z,$\textit{\ }$f($\textit{\textperiodcentered }$,x)$%
\textit{\ is measurable on }$U^{x}$\textit{\ and for each }$\omega \in
\Omega ,$\textit{\ }$f(\omega ,$\textit{\textperiodcentered }$)$\textit{\ is
continuous on }$U^{\omega }.$\textit{\ Moreover, }$f($\textit{%
\textperiodcentered }$,$\textit{\textperiodcentered }$)$\textit{\ is jointly
measurable.\medskip }

\begin{theorem}
(U. s. c. Lifting Theorem. [30]\textbf{)}. \textit{Let }$Y$\textit{\ be a
separable space, }$(\Omega $\textit{, }$\tciFourier ,\mu )$\textit{\ be a
complete finite measure space and }$X:\Omega \rightarrow 2^{Y}$\textit{\ be
an integrably bounded, nonempty, convex valued correspondence such that for
all }$\omega \in \Omega ,$\textit{\ }$X(\omega )$\textit{\ is a weakly
compact, convex subset of }$Y$\textit{. Denote by }$S_{X}$\textit{\ the set }%
$\{x\in L_{1}(\mu ,Y):x(\omega )\in X(\omega )$\textit{\ }$\mu -a.e.\}.$%
\textit{\ Let }$T:\Omega \times S_{X}\rightarrow 2^{Y}$\textit{\ be a
nonempty, closed, convex valued correspondence such that }$T(\omega
,x)\subset X(\omega )$\textit{\ for all }$(\omega ,x)\in \Omega \times
S_{X}^{1}.$\textit{\ Assume that for each fixed }$x\in S_{X},$\textit{\ }$T($%
\textit{\textperiodcentered }$,x)$\textit{\ has a measurable graph and that
for each fixed }$\omega \in \Omega ,$\textit{\ }$T(\omega ,$\textit{%
\textperiodcentered }$):S_{X}\rightarrow 2^{Y}$\textit{\ is u.s.c. in the
sense that the set }$\{x\in S_{X}:T(\omega ,x)\subset V\}$\textit{\ is
weakly open in }$S_{X}$\textit{\ for every norm open subset }$V$\textit{\ of 
}$Y$\textit{. Define the correspondence }$\Phi :S_{X}\rightarrow 2^{S_{X}}$%
\textit{\ by}
\end{theorem}

$\Phi (x)=\{y\in S_{X}:y(\omega )\in T(\omega ,x)$\textit{\ }$\mu -a.e.\}.$

\textit{Then }$\Phi $\textit{\ is weakly u.s.c., i.e., the set }$\{x\in
S_{X}:\Phi (x)\subset V\}$\textit{\ is weakly open in }$S_{X}$\textit{\ for
every weakly open subset }$V$\textit{\ of }$S_{X}$\textit{.}

\begin{center}
\bigskip
\end{center}

\end{document}